\pdfoutput=1
\documentclass[twoside]{article}
\usepackage[letterpaper, left=4.7cm, right=4.7cm, top=4cm, bottom=3.5cm]{geometry}
\usepackage{graphicx}
\usepackage{amsmath,amssymb,amsthm}
\usepackage{authblk}
\usepackage[utf8]{inputenc}
\usepackage{microtype}
\usepackage{mathrsfs}  
\usepackage{enumitem}
\usepackage{calc}
\usepackage{esint}
\usepackage{fancyhdr}
\usepackage{xcolor}
\usepackage[labelfont = bf]{caption}
\usepackage{doi}
\usepackage{subfigure}
\usepackage[numbers]{natbib}
\usepackage{physics}
\usepackage{tikz}
\usepackage{mathdots}
\usepackage{yhmath}
\usepackage{cancel}
\usepackage{color}
\usepackage{siunitx}
\usepackage{array}
\usepackage{multirow}
\usepackage{gensymb}
\usepackage{tabularx}
\usepackage{extarrows}
\usepackage{booktabs}
\usetikzlibrary{fadings}
\usetikzlibrary{patterns}
\usetikzlibrary{shadows.blur}
\usetikzlibrary{shapes}

\newtheorem{theorem}{Theorem}[section]
\numberwithin{equation}{section}
\newtheorem{proposition}[theorem]{Proposition}

\newtheorem{corollary}[theorem]{Corollary}
\newtheorem{remark}[theorem]{Remark}
\newtheorem{lemma}[theorem]{Lemma}
\newtheorem{example}[theorem]{Example}
\newtheorem{algorithm}[theorem]{Algorithm}

\usepackage{titlesec}
\titleformat{\section}{\normalfont\scshape\centering}{\thesection.}{0.5em}{}
\titleformat*{\subsection}{\itshape}

\providecommand{\keywords}[1]
{
	{\small\textit{Keywords:~~} #1}
}

\providecommand{\MSC}[1]
{
	{\small\textit{AMS MSC (2010):~~} #1}
}

\usepackage{hyperref}
\hypersetup{
	colorlinks=true,
	linkcolor=blue,
	citecolor = blue,
	filecolor=magenta,      
	urlcolor=cyan,
	pdfpagemode=FullScreen,
}

\makeatletter
\renewcommand{\NAT@separator}{\NAT@sep\nolinebreak}
\makeatother

\begin{document}
	\setlength{\abovedisplayskip}{5.5pt}
	\setlength{\belowdisplayskip}{5.5pt}
	\setlength{\abovedisplayshortskip}{5.5pt}
	\setlength{\belowdisplayshortskip}{5.5pt}

	\title{Error \hspace*{-0.5mm}estimates \hspace*{-0.5mm}for total-variation \hspace*{-0.5mm}regularized \mbox{minimization \hspace*{-0.5mm}problems \hspace*{-0.5mm}with \hspace*{-0.5mm}singular \hspace*{-0.5mm}dual \hspace*{-0.5mm}solutions}}
	\author[1]{Sören Bartels\thanks{Email: \texttt{bartels@mathematik.uni-freiburg.de}}}
	\author[2]{Alex Kaltenbach\thanks{Email: \texttt{alex.kaltenbach@mathematik.uni-freiburg.de}}}
	\date{\today}
	\affil[1]{\small{Institute of Applied Mathematics, Albert--Ludwigs--University Freiburg, Hermann--Herder--Stra\ss e 10, 79104 Freiburg}}
		\affil[2]{\small{Institute of Applied Mathematics, Albert--Ludwigs--University Freiburg, Ernst--Zermelo--Stra\ss e 1, 79104 Freiburg}}
	\maketitle

	\pagestyle{fancy}
	\fancyhf{}
	\fancyheadoffset{0cm}
	\addtolength{\headheight}{-0.25cm}
	\renewcommand{\headrulewidth}{0pt} 
	\renewcommand{\footrulewidth}{0pt}
	\fancyhead[CO]{\textsc{Error estimates for total-variation regularized minimization}}
	\fancyhead[CE]{\textsc{S. Bartels and A. Kaltenbach}}
	\fancyhead[R]{\thepage}
	\fancyfoot[R]{}
	
	\begin{abstract}
		\hspace{-0.25em}Recent quasi-optimal error estimates for the finite element \mbox{approximation} of total-variation regularized minimization problems using the Crouzeix--Raviart finite element require the existence of a Lipschitz continuous~dual solution, which 
		is not generally given. We provide  analytic proofs showing that the Lipschitz continuity of a dual solution~is not~necessary,~in~\mbox{general}. 
		Using the Lipschitz truncation technique, we, in  addition, derive~error~estimates that depend directly on the Sobolev regularity~of~a~given~dual~\mbox{solution}.~~~~~~~~~~~  
	\end{abstract}

	
	\MSC{26A45, 65N15, 65N30, 68U10.}
	
	\section{Introduction}\label{sec:intro}
		\thispagestyle{empty}
	\qquad In this article, we examine the finite element discretization of the Rudin--Osher--Fatemi (ROF) model from \cite{ROF92}, which serves as a model problem for general convex and non-smooth minimization problems. This image processing model  determines a function $u\in BV(\Omega)\cap L^2(\Omega)$ via minimizing $I:BV(\Omega)\cap L^2(\Omega)\to \mathbb{R}$, defined by 
	\begin{align}
		I(u):=\vert \textup{D}u\vert (\Omega)+\frac{\alpha}{2}\|u-g\|^2_{L^2(\Omega)}\label{intro-primal}
	\end{align}
	for all $u\in BV(\Omega)\cap L^2(\Omega)$, where $\vert \textup{D}u\vert (\Omega)$ denotes the total variation, $g\in L^2(\Omega)$ is the input data, e.g., a noisy image, and $\|u-g\|_{L^2(\Omega)}^2$ is the so-called fidelity~term. In addition, the fidelity parameter $\alpha>0$ is a given constant, which determines the balance between de-noising and preserving the input image. For a more in-depth analysis of this model, concerning its
 	analytical properties,~explicit~solutions, and numerical~\hspace*{-0.1mm}methods,~\hspace*{-0.1mm}we~\hspace*{-0.1mm}\mbox{refer}~\hspace*{-0.1mm}to~\hspace*{-0.1mm}{\cite{CL97,AFP00,HK04,CCN07,Sparse10,CLL11,WL11,Bar12,LM12,ABM14,BNS15,HHSVW19,CP20,Bar21,BTW21}}.
	Since this model allows for and preserves discontinuities of the input~data~$g$,~cf.~\cite{CCN07}, continuous finite element methods are known to perform sub-optimally,~cf.~\cite{BNS15,Bar21}. Recent contributions, cf. \cite{CP20,Bar20,Bar21}, reveal that~the~\mbox{quasi-optimal}~\mbox{convergence}~rate $\smash{\mathcal{O}(h^{\frac{1}{2}})}$ for discontinuous solutions on quasi-uniform triangulations can be obtained using discontinuous, low-order Crouzeix--Raviart finite~elements~introduced~in~\cite{CR73}. More precisely, these error estimates yield a bound for the error~for~the~approximation of
	minimizers of $I\hspace*{-0.1em}:\hspace*{-0.1em}BV(\Omega)\cap L^2(\Omega)\hspace*{-0.1em}\to\hspace*{-0.1em} \mathbb{R}$ via minimizing~the~discrete~functional $I_h:\mathcal{S}^{1,cr}(\mathcal{T}_h)\to \mathbb{R}$, defined by
	\begin{align*}
		I_h(u_h):=\|\nabla_h u_h\|_{L^1(\Omega;\mathbb{R}^d)}+\frac{\alpha}{2}\|\Pi_h(u_h-g)\|_{L^2(\Omega)}^2
	\end{align*}
	for all $u_h\hspace*{-0.1em}\in \hspace*{-0.1em}\mathcal{S}^{1,cr}(\mathcal{T}_h)$, where $\mathcal{S}^{1,cr}(\mathcal{T}_h)$ is the Crouzeix--Raviart finite~element~space, i.e., the space of piece-wise affine functions that are continuous at the midpoints of element sides,  $\nabla_h\!:\!\mathcal{S}^{1,cr}(\mathcal{T}_h)\!\to\! \mathcal{L}^0(\mathcal{T}_h)^d$ denotes the element-wise~\mbox{gradient},~and~$\Pi_h\!:\!L^2(\Omega)\!\hspace*{-0.2mm}\to\!\hspace*{-0.2mm} \mathcal{L}^0(\mathcal{T}_h)$ \hspace*{-0.1mm}is \hspace*{-0.1mm}the \hspace*{-0.1mm}$L^2$\hspace*{-0.2mm}--projection~\hspace*{-0.1mm}operator~\hspace*{-0.1mm}onto~\hspace*{-0.1mm}\mbox{element-wise}~\hspace*{-0.1mm}constant~\hspace*{-0.1mm}\mbox{functions}. Note that the family of discrete functionals $I_h:\mathcal{S}^{1,cr}(\mathcal{T}_h)\to \mathbb{R}$, $h>0$, defines a non-conforming approximation of the functional $I:BV(\Omega)\cap L^2(\Omega)\to \mathbb{R}$,~as,~e.g., jump \hspace*{-0.1mm}terms \hspace*{-0.1mm}of \hspace*{-0.1mm}$u_h$ \hspace*{-0.2mm}across \hspace*{-0.1mm}inter-element \hspace*{-0.1mm}sides \hspace*{-0.1mm}are \hspace*{-0.1mm}not \hspace*{-0.1mm}included.~\hspace*{-0.2mm}For~\hspace*{-0.1mm}this~\hspace*{-0.1mm}\mbox{family}~\hspace*{-0.1mm}\mbox{recently} a \hspace*{-0.1mm}$\Gamma\!$--convergence \hspace*{-0.1mm}result  \hspace*{-0.1mm}with \hspace*{-0.1mm}respect \hspace*{-0.1mm}to \hspace*{-0.1mm}strong \hspace*{-0.1mm}convergence~\hspace*{-0.1mm}in~\hspace*{-0.1mm}$L^1(\Omega)$~\hspace*{-0.1mm}or~\hspace*{-0.1mm}\mbox{distributional} convergence  has been established under general assumptions, i.e., that $g\in L^2(\Omega)$, cf. \hspace*{-0.3em}\cite[\hspace*{-0.05em}Propositon~\hspace*{-0.05em}3.1]{CP20}. \!\hspace*{-0.1em}However, the quasi-optimal~rate~$\smash{\mathcal{O}(h^{\frac{1}{2}}\hspace*{-0.1em})}$~till~now~only~holds~if the dual problem~given~via~maximizing $D:W^2_N(\textup{div};\Omega)\cap L^\infty(\Omega;\mathbb{R}^d)\to \mathbb{R}\cup\{-\infty\}$, defined~by 
	\begin{align}
		D(z):=-\frac{1}{2\alpha}\|\textup{div}(z)+\alpha g\|_{L^2(\Omega)}^2+\frac{\alpha}{2}\|g\|_{L^2(\Omega)}^2-I_{K_1(0)}(z)\label{intro-dual}
	\end{align}
	for all $z\in W^2_N(\textup{div};\Omega)\cap L^\infty(\Omega;\mathbb{R}^d)$, where ${I_{K_1(0)}:L^\infty(\Omega;\mathbb{R}^d)\to \mathbb{R}\cup\{+\infty\}}$~is for $z\hspace*{-0.1em}\in \hspace*{-0.17em} L^\infty(\Omega;\mathbb{R}^d)$ defined by $\hspace*{-0.1em}I_{K_1(0)}(z)\hspace*{-0.17em}:=\hspace*{-0.17em}0$ if $\|z\|_{L^\infty(\Omega;\mathbb{R}^d)}\hspace*{-0.17em}\leq\hspace*{-0.17em} 1$~and~${I_{K_1(0)}(z)\hspace*{-0.17em}:=\hspace*{-0.17em}+\infty}$~else, admits a Lipschitz continuous solution. Unfortunately, the Lipschitz continuity of a maximum of $D:W^2_N(\textup{div};\Omega)\cap L^\infty(\Omega;\mathbb{R}^d)\to \mathbb{R}\cup\{-\infty\}$ is not generally~given,~as 
	\cite[Section 3]{BTW21} clarified. Without  imposing  the existence~of~a~Lipschitz~continuous solution to \eqref{intro-dual}, but  that $g\!\in\! BV(\Omega)\cap L^\infty(\Omega)$, in \cite[Section~5.2]{CP20},~the~\mbox{sub-optimal} convergence rate $\smash{\mathcal{O}(h^{\frac{1}{4}})}$
	has been established. The approach of \cite[Section~5.2]{CP20} consists in a convolution
	of a maximum $z\in W^\infty_N(\textup{div};\Omega)$~of~\eqref{intro-dual}~in~order~to~comply with the crucial Lipschitz continuity property at least~in an approximate~sense.\linebreak 
	 We \hspace*{-0.1mm}use \hspace*{-0.1mm}an \hspace*{-0.1mm}alternative \hspace*{-0.1mm}regularization \hspace*{-0.1mm}approach, \hspace*{-0.2mm}which \hspace*{-0.1mm}operates \hspace*{-0.1mm}highly~\hspace*{-0.1mm}at~\hspace*{-0.1mm}a~\hspace*{-0.1mm}\mbox{local}~\hspace*{-0.1mm}level, the celebrated Lipschitz truncation technique.~Its~basic~purpose is to approximate  Sobolev functions $u \in  W^{1,p}(\Omega)$~by~$\lambda$--Lipschitz~functions~${u_\lambda\in W^{1,\infty}(\Omega)}$,~${\lambda>0}$.
	 The original approach of this technique traces back~to~Acerbi~and~Fusco,~cf.~\cite{AF84,AF87,AF88}. Since then, the Lipschitz truncation technique is used in various areas of analysis: In \hspace*{-0.1em}the \hspace*{-0.1em}calculus \hspace*{-0.1em}of \hspace*{-0.1em}variations, \hspace*{-0.1em}in \hspace*{-0.1em}the \hspace*{-0.1em}existence~\hspace*{-0.1em}theory~\hspace*{-0.1em}of~\hspace*{-0.1em}partial~\hspace*{-0.1em}differential~\hspace*{-0.1em}\mbox{equations}, and in regularity theory. For a longer list of references,~we~refer~the~reader~to~\cite{DMS08}. To the best of the authors knowledge, this article provides the first deployment of the Lipschitz truncation technique in the field of image processing.~More~\mbox{precisely}, for the application in this article,
	 the main advantage~of~the~Lipschitz~\mbox{truncation} technique in comparison to convolution is that not only $u$ and $u_\lambda$ coincide up to a set of small measure, but equally~$\nabla u$~and~$\nabla u_\lambda$~do.~By~deploying~the~Lipschitz~truncation technique, we arrive at error estimates whose resulting rates directly depend on the respective Sobolev regularity of a given maximum~${z\!\in\! W^{1,p}(\Omega;\mathbb{R}^d)}$~of~\eqref{intro-dual}. 
	 If only $g\in L^\infty(\Omega)$ and one additionally has that, e.g.,
	 ${z\in W^{1,p}(\Omega;\mathbb{R}^d)}$~for~${p\ge 3}$, then the results of this article~yield~the~\mbox{sub-optimal}~rate~$\smash{\mathcal{O}(h^\frac{1}{4})}$. 
	 In this manner, \linebreak we intend to fill the gap between the optimal rate $\smash{\mathcal{O}(h^{\frac{1}{2}})}$ for $z\in W^{1,\infty}(\Omega;\mathbb{R}^d)$ and $g\in L^\infty(\Omega)$ and the  rate $\smash{\mathcal{O}(h^{\frac{1}{4}})}$ for $z\in W^\infty_N(\textup{div};\Omega)$ and~${g\in L^\infty(\Omega)\cap BV(\Omega)}$.~~~~~~~~
	  
	As a maximum of \eqref{intro-dual} is not necessarily in a Sobolev~space,~but~in~$W^2_N(\textup{div};\Omega)$ $\cap L^\infty(\Omega;\mathbb{R}^d)$, we also study the case of a non-existence of Sobolev~solutions~to~\eqref{intro-dual}. It~turns out that if a maximum $z\!\in\!  W^2_N(\textup{div};\Omega)\cap L^\infty(\Omega;\mathbb{R}^d)$~of~\eqref{intro-dual}~is~\mbox{element-wise} Lipschitz continuous, i.e., the discontinuity set $J_z$ is resolved~by~the~\mbox{triangulations}, or at least in an approximate sense with the rate $\mathcal{O}(h)$, cf. Remark~\ref{rem:soeren_trick_special_cases}, then the optimal rate $\smash{\mathcal{O}(h^{\frac{1}{2}})}$ can be expected. Beyond that, we find that the~optimal~rate $\smash{\mathcal{O}(h^{\frac{1}{2}})}$ is attained if a dual solution fulfills $\vert z\vert<1$ along its discontinuity set $J_z$ while, simultaneously, its jump $[\![z]\!]$ over its discontinuity set  $J_z$ remains~small.
	Some of these conditions apply, e.g., to the setting described in \cite[Section 3]{BTW21} with a suitable triangulation $\mathcal{T}_h$, $h>0$, of the domain $\Omega$, for which the optimal rate $\smash{\mathcal{O}(h^{\frac{1}{2}})}$ could be reported~without~giving~an~analytical~explanation. This article's purpose is to give -- at least for special cases -- a~missing~analytical~explanation.~~~~~

	\textit{This article is organized as follows:} In Section \ref{sec:preliminaries}, we introduce the employed notation, define the relevant finite element spaces and give a brief review of the continuous and discretized ROF model. In Section \ref{sec:lipschitz_truncation},  using the Lipschitz~\mbox{truncation} technique, we establish error estimates that depend directly on the Sobolev regularity of a maximum of \eqref{intro-dual}. In Section \ref{sec:non_sobolev}, we prove quasi-optimal error estimates without explicitly imposing that  a Lipschitz continuous~maximum~of~\eqref{intro-dual}~exists. In Section \ref{sec:experiments}, we confirm our theoretical findings~via~numerical~experiments.
	
	\section{Preliminaries}\label{sec:preliminaries}
	
	\qquad Throughout the article, if not otherwise specified, we denote by ${\Omega\hspace*{-0.1em}\subseteq\hspace*{-0.1em} \mathbb{R}^d}$,~${d\hspace*{-0.1em}\in\hspace*{-0.1em}\mathbb{N}}$, a bounded polyhedral Lipschitz domain, whose boundary is disjointly divided~into a Dirichlet part $\Gamma_D$ and a Neumann part $\Gamma_N$, i.e.,  ${\partial\Omega\hspace*{-0.1em}=\hspace*{-0.1em}\Gamma_D\cup\Gamma_N}$~and~${\emptyset\hspace*{-0.1em}=\hspace*{-0.1em}\Gamma_D\cap\Gamma_N}$.
	\subsection{Function spaces}\label{subsec:function_spaces}

	\qquad For $p\in \left[1,\infty\right]$ and $l\in \mathbb{N}$, we employ the standard notations~~~~~~~~~~~~~~~~~~~~~~~~~~~~~
	\begin{align*}
		W^{1,p}_D(\Omega;\mathbb{R}^l)&:=\big\{u\in L^p(\Omega;\mathbb{R}^l)\mid \nabla u\in L^p(\Omega;\mathbb{R}^{l\times d}),\, \textup{tr}(u)=0\text{ in }L^p(\Gamma_D;\mathbb{R}^l)\big\},\\
	W^{p}_N(\textup{div};\Omega)&:=\big\{z\in L^p(\Omega;\mathbb{R}^d)\mid \textup{div}(z)\in L^p(\Omega),\,\textup{tr}(z)\cdot n=0\text{ in }W^{-\frac{1}{p},p}(\Gamma_N)\big\},
	\end{align*}
	$\smash{W^{1,p}(\Omega;\mathbb{R}^l):=W^{1,p}_D(\Omega;\mathbb{R}^l)}$ if $\smash{\Gamma_D=\emptyset}$, and $\smash{W^{p}(\textup{div};\Omega):=W^{p}_N(\textup{div};\Omega)}$ if $\smash{\Gamma_N=\emptyset}$,
	where $\smash{\textup{tr}\hspace*{-0.05em}:\hspace*{-0.05em}W^{1,p}(\Omega;\mathbb{R}^l)\hspace*{-0.05em}\to\hspace*{-0.05em} L^p(\partial\Omega)}$  and $
	\smash{\textup{tr}(\cdot)\cdot n\hspace*{-0.05em}:\hspace*{-0.05em}W^p(\textup{div};\Omega)\hspace*{-0.05em}\to\hspace*{-0.05em} (W^{1,p'}(\Omega))^*}$ 
	denote the trace and normal trace operator. In particular, we  \mbox{predominantly}~\mbox{omit} $\textup{tr}(\cdot)$ in this context. Apart from that, we fall back on the abbreviations~${L^p(\Omega) \!:= \!L^p\hspace*{-0.05em}(\Omega;\hspace*{-0.05em}\mathbb{R}^1\hspace*{-0.05em})}$ and  $W^{1,p}(\Omega):=W^{1,p}(\Omega;\mathbb{R}^1)$. Let  $\vert \textup{D}(\cdot)\vert(\Omega):L^1_{\textup{loc}}(\Omega)\to \mathbb{R}\cup\{+\infty\}$, defined by\footnote{Here, $C_c^\infty(\Omega;\mathbb{R}^d)$ denotes the space of smooth and in $\Omega$ compactly supported vector  fields.\vspace*{-5mm}} 
	\begin{align*}
		\vert \textup{D}u\vert(\Omega):=\sup_{\phi\in C_c^\infty(\Omega;\mathbb{R}^d),\|\phi\|_{L^\infty(\Omega;\mathbb{R}^d)\leq 1}}{-\int_{\Omega}{u\,\textup{div}(\phi)\,\textrm{d}x}}
	\end{align*}
	for all $u\in L^1_{\textup{loc}}(\Omega)$, denote the total variation. Then, the space of functions of bounded variation is defined by $BV(\Omega):=\big\{u\in L^1(\Omega)\mid \vert \textup{D}u\vert(\Omega)<\infty\big\}$.
	\begin{align*}
		BV(\Omega):=\big\{u\in L^1(\Omega)\mid \vert \textup{D}u\vert(\Omega)<\infty\big\}.
	\end{align*}

	\subsection{Triangulations}
	
	\qquad In what follows, we let $(\mathcal{T}_h)_{h>0}$ be a sequence of  regular, i.e., uniformly shape regular and conforming, triangulations of $\Omega\subseteq \mathbb{R}^d$, $d\in\mathbb{N}$, cf. \cite{BS08}. The sets $\mathcal{S}_h$ and $\mathcal{N}_h$  contain the sides and vertices, resp., of the elements. The parameter $h>0$ refers to  the maximal mesh-size of $\mathcal{T}_h$. More precisely,~if~we~define~${h_T:=\textup{diam}(T)}$  for all ${T\in \mathcal{T}_h}$, then we have that $h=\max_{T\in \mathcal{T}_h}{h_T}$. For any  $k\in \mathbb{N}$~and~${T\in \mathcal{T}_h}$, we let $\mathcal{P}_k(T)$ denote the set of polynomials of maximal total degree $k$ on $T$. Then, the set of element-wise polynomial functions or vector fields, resp., is defined by
	\begin{align*}
		\smash{\mathcal{L}^k(\mathcal{T}_h)^l:=\big\{v_h\in L^\infty(\Omega;\mathbb{R}^l)\mid v_h|_T\in\mathcal{P}_k(T)\text{ for all }T\in \mathcal{T}_h\big\}.}
	\end{align*}
	For any $T\in \mathcal{T}_h$ and $S\in \mathcal{S}_h$, we let   $\smash{x_T:=\frac{1}{d+1}\sum_{z\in \mathcal{N}_h\cap T}{z}}$ and $\smash{x_S:=\frac{1}{d}\sum_{z\in \mathcal{N}_h\cap S}{z}}$ denote the midpoints (barycenters) of $T$ and $S$, resp. The $L^2$--projection~operator onto piece-wise constant functions or vector fields, resp., is denoted by
	\begin{align*}
		\smash{\Pi_h:L^1(\Omega;\mathbb{R}^l)\to \mathcal{L}^0(\mathcal{T}_h)^l.}
	\end{align*}
	For $v_h\hspace*{-0.1em}\in \hspace*{-0.1em} \mathcal{L}^1(\mathcal{T}_h)^l$, it holds $\Pi_hv_h|_T\hspace*{-0.1em}=\hspace*{-0.1em}v_h(x_T)$ for all $T\hspace*{-0.1em}\in\hspace*{-0.1em} \mathcal{T}_h$. Moreover,~for~${p\hspace*{-0.1em}\in\hspace*{-0.1em} \left[1,\infty\right]}$, there exists a constant $c_{\Pi}>0$ such that for all $v\in L^p(\Omega;\mathbb{R}^l)$, cf. \cite{EG04}, we~have~that
	\begin{description}[noitemsep,topsep=1.5pt,font=\normalfont\itshape]
		\item[(L0.1)]\hypertarget{(L0.1)}{} $\|\Pi_h v\|_{L^p(\Omega;\mathbb{R}^l)}\leq \| v\|_{L^p(\Omega;\mathbb{R}^l)}$,
		\item[(L0.2)]\hypertarget{(L0.2)}{} $\|v-\Pi_h v\|_{L^p(\Omega;\mathbb{R}^l)}\leq c_{\Pi}h \|\nabla v\|_{L^p(\Omega;\mathbb{R}^{l\times d})}$ if $v\in W^{1,p}(\Omega;\mathbb{R}^l)$.
	\end{description}
	
	\subsection{Crouzeix--Raviart finite elements}
	
	\qquad A particular instance of a larger class of non-conforming finite element~spaces, introduced in \cite{CR73},
	is the Crouzeix--Raviart finite element space, which~\mbox{consists}~of piece-wise affine functions that are continuous at the midpoints of~element~sides,~i.e.,
	\begin{align*}
		\smash{\mathcal{S}^{1,cr}(\mathcal{T}_h):=\big\{v_h\in \mathcal{L}^1(\mathcal{T}_h)\mid v_h\text{ is continuous in }x_S\text{ for all }S\in \mathcal{S}_h\big\}.}
	\end{align*}
	The element-wise application of the gradient to $v_h\in \mathcal{S}^{1,cr}(\mathcal{T}_h)$ defines an element-wise constant vector field $\nabla_hv_h\in \mathcal{L}^0(\mathcal{T}_h)^d$ via ${\nabla_hv_h|_T:=\nabla(v_h|_T)}$~for~all~${T\in \mathcal{T}_h}$. Crouzeix--Raviart finite element functions that vanish at midpoints of boundary~element \hspace*{-0.1mm}sides \hspace*{-0.1mm}that \hspace*{-0.1mm}correspond \hspace*{-0.1mm}to \hspace*{-0.1mm}the \hspace*{-0.1mm}Dirichlet \hspace*{-0.1mm}boundary \hspace*{-0.1mm}$\Gamma_{\! D}$ \hspace*{-0.1mm}are \hspace*{-0.1mm}contained~\hspace*{-0.1mm}in~\hspace*{-0.1mm}the~\hspace*{-0.1mm}space
	\begin{align*}
	\smash{	\mathcal{S}^{1,cr}_D(\mathcal{T}_h):=\big\{v_h\in\mathcal{S}^{1,cr}(\mathcal{T}_h)\mid v_h(x_S)=0\text{ for all }S\in \mathcal{S}_h\text{ with }S\subseteq \Gamma_D\big\}.}
	\end{align*}
	In particular, we have that $	\mathcal{S}^{1,cr}_D(\mathcal{T}_h)=	\mathcal{S}^{1,cr}(\mathcal{T}_h)$ if $\Gamma_D=\emptyset$.
	A basis of  $\mathcal{S}^{1,cr}(\mathcal{T}_h)$ is given by the functions $\varphi_S\in \mathcal{S}^{1,cr}(\mathcal{T}_h)$, $S\in \mathcal{S}_h$, satisfying the Kronecker~property $\varphi_S(x_{S'})=\delta_{S,S'}$ for all $S,S'\in \mathcal{S}_h$. A basis of  $\mathcal{S}^{1,cr}_D(\mathcal{T}_h)$~is~given~by~$(\varphi_S)_{S\in \mathcal{S}_h;S\not\subseteq\Gamma_D}$. For any $p\in \left[1,\infty\right]$, the quasi-interpolation operator ${I_{cr}:W^{1,p}_D(\Omega)\to \mathcal{S}^{1,cr}_D(\mathcal{T}_h)}$, for all $v\in W^{1,p}_D(\Omega)$ defined by\vspace*{-2mm}
	\begin{align}
		I_{cr}v:=\smash{\sum_{S\in \mathcal{S}_h}{v_S\varphi_S},}\quad v_S:=\fint_S{v\,\textup{d}s}\label{CR-interpolant}
	\end{align}
	preserves  averages of gradients, i.e., $\nabla_h(I_{cr}v)\!=\!\Pi_h(\nabla v)$ in  $\mathcal{L}^0(\mathcal{T}_h)^d$  for ${v\!\in\! W^{1,p}_D(\Omega)}$. Moreover, for $p\hspace*{-0.15em}\in\hspace*{-0.15em} \left[1,\infty\right]$, there exits a constant ${c_{cr}\hspace*{-0.15em}>\hspace*{-0.15em}0}$~such~that for all ${v\hspace*{-0.15em}\in\hspace*{-0.15em} W^{1,p}_D(\Omega)}$, cf. \cite{BBF13}, we have that
	\begin{description}[noitemsep,topsep=1.5pt,font=\normalfont\itshape]
		\item[(CR.1)]\hypertarget{(CR.1)}{} $\|\nabla_h (I_{cr}v)\|_{L^p(\Omega;\mathbb{R}^d)}\leq \|\nabla v\|_{L^p(\Omega;\mathbb{R}^d)}$,
		\item[(CR.2)]\hypertarget{(CR.2)}{} $\|v-I_{cr}v\|_{L^p(\Omega)}\leq c_{cr}h \|\nabla v\|_{L^p(\Omega;\mathbb{R}^d)}$,
		\item[(CR.3)]\hypertarget{(CR.3)}{} $\|I_{cr}v\|_{L^\infty(\Omega)}\leq c_d\| v\|_{L^\infty(\Omega)}$, where $c_d:=(d+1)(d-1)$,   if $v\in L^\infty(\Omega)$.
	\end{description}
	For \hspace*{-0.1mm}$p\!=\!1$, \hspace*{-0.1mm}due \hspace*{-0.1mm}to \hspace*{-0.1mm}the \hspace*{-0.1mm}density \hspace*{-0.1mm}of \hspace*{-0.1mm}$C^\infty(\Omega)\cap BV(\Omega)$ \hspace*{-0.1mm}in \hspace*{-0.1mm}$BV(\Omega)$,~\hspace*{-0.1mm}cf.~\hspace*{-0.1mm}\cite{ABM14},~\hspace*{-0.1mm}the~\hspace*{-0.1mm}\mbox{operator}~\hspace*{-0.1mm}and \textit{(\hyperlink{(CR.1)}{CR.1})}--\textit{(\hyperlink{(CR.3)}{CR.3})} can be extended to $v\in BV(\Omega)$, losing~the~representation~\eqref{CR-interpolant}.\vspace*{-5mm}
	\newpage
	
	\subsection{Raviart--Thomas finite elements}

	\qquad The lowest order Raviart--Thomas finite element space, introduced in \cite{RT75}, consists of \mbox{piece-wise} affine vector fields that possess weak divergences, i.e.,
	\begin{align*}
		\smash{\mathcal{R}T^0(\mathcal{T}_h):=\big\{z_h\in \mathcal{L}^1(\mathcal{T}_h)^d\mid z_h|_T\cdot n_T=-z_h|_{T'}\cdot n_{T'}\text{ on }T\cap T'\text{ if }T\cap T\in \mathcal{S}_h\big\},}
	\end{align*}
	where $n_T:\partial T\to \mathbb{S}^{d-1}$ for all $T\in \mathcal{T}_h$ denotes the unit normal vector field to $ T$ pointing outward. Raviart--Thomas finite element~functions~that~have~vanishing normal components on the Neumann boundary $\Gamma_N$ are~contained~in~the~space 
	\begin{align*}
		\smash{\mathcal{R}T^0_N(\mathcal{T}_h):=\big\{z_h\in	\mathcal{R}T^0(\mathcal{T}_h)\mid z_h\cdot n=0\text{ on }\Gamma_N\big\}.}
	\end{align*}
	In particular, we have that $\mathcal{R}T^0_N(\mathcal{T}_h)=\mathcal{R}T^0(\mathcal{T}_h)$ if $\Gamma_N=\emptyset$. A basis of  $\mathcal{R}T^0(\mathcal{T}_h)$ is given by the vector fields $\psi_S\in \mathcal{R}T^0(\mathcal{T}_h)$, $S\in \mathcal{S}_h$, satisfying the Kronecker property $\psi_S|_{S'}\cdot n_{S'}=\delta_{S,S'}$ on $S'$ for all $S\in \mathcal{S}_h$, where $n_S$ for all $S\in \mathcal{S}_h$ denotes the unit normal vector on $S$ that points from $T_-$ to $T_+$ if $S\!=\!\partial T_-\cap \partial T_+\!\in\! \mathcal{S}_h$.~A~basis~of  $\mathcal{R}T^0_N(\mathcal{T}_h)$ is given by 
	 $\psi_S\!\in\! \mathcal{R}T^0_N(\mathcal{T}_h)$,~${S\!\in\! \mathcal{S}_h\!\setminus\!\Gamma_N}$.~The~\mbox{quasi-interpolation}~\mbox{operator} $I_{\mathcal{R}T}:V^{\textup{div}}(\Omega):=\{y\in L^p(\Omega;\mathbb{R}^d) \mid \textup{div}(y)\in L^q(\Omega)\}\to  \mathcal{R}T^0_N(\mathcal{T}_h)$, where $p>2$ and $\smash{q>\frac{2d}{d+2}}$, for~all~$\smash{z\in V^{\textup{div}}(\Omega)}$ defined~by\vspace*{-2.5mm}
	\begin{align}
		I_{\mathcal{R}T}z:=\smash{\sum_{S\in \mathcal{S}_h}{z_S\psi_S}},\quad z_S:=\fint_S{z\cdot n_S\,\textup{d}s}\label{RT-interpolant}
	\end{align}
	preserves averages of divergences, i.e.,  $\textup{div}(I_{\mathcal{R}T}z)=\Pi_h(\textup{div}(z))$ in $\mathcal{L}^0(\mathcal{T}_h)$ for all ${z\in V^{\textup{div}}(\Omega)}$. Moreover, for $p\in \left[1,\infty\right]$, there exists~a~constant ${c_{\mathcal{R}T}>0}$~such~that for all $z\in V^{\textup{div}}(\Omega)$,~cf.~\cite{EG21}, we have that
	\begin{description}[noitemsep,topsep=0.5pt,font=\normalfont\itshape]
		\item[(RT.1)]\hypertarget{(RT.1)}{} $\|z-I_{\mathcal{R}T}z\|_{L^p(\Omega;\mathbb{R}^d)}\leq c_{\mathcal{R}T}h \|\nabla z\|_{L^p(\Omega;\mathbb{R}^{d\times d})}$ if $z\in W^{1,p}(\Omega;\mathbb{R}^d)$,
		\item[(RT.2)]\hypertarget{(RT.2)}{} $\|I_{\mathcal{R}T}z\|_{L^\infty(\Omega;\mathbb{R}^d)}\leq c_{\mathcal{R}T} \| z\|_{L^\infty(\Omega;\mathbb{R}^d)}$ if $z\in L^\infty(\Omega;\mathbb{R}^d)$.
	\end{description}
	For \hspace*{-0.1mm}$p\!\in\! [1,\infty)$, \hspace*{-0.2mm}due \hspace*{-0.1mm}to
	\hspace*{-0.1mm}the \hspace*{-0.1mm}density \hspace*{-0.1mm}of \hspace*{-0.1mm}$C^\infty(\Omega;\mathbb{R}^d)\cap W^p_N(\textup{div};\Omega)$ \hspace*{-0.1mm}in  \hspace*{-0.2mm}$W^p_N(\textup{div};\Omega)$,~\hspace*{-0.1mm}the~\hspace*{-0.1mm}oper-ator and \textit{(\hyperlink{(RT.2)}{RT.2})} can be extended to $z\!\in\! W^p_N(\textup{div};\Omega)$, losing the~\mbox{representation}~\eqref{RT-interpolant}.

	\subsection{The continuous  Rudin--Osher--Fatemi (ROF) model}\label{subsec:continuous_ROF_model}
	
	\qquad Given $g\!\in\! L^2(\Omega)$ and $\alpha\!>\!0$, the Rudin--Osher--Fatemi~(ROF)~model,~cf.~\cite{ROF92},~determines a function ${u\!\in\! BV(\Omega)\cap L^2(\Omega)}$ that is minimal for ${I\!:\!BV(\Omega)\cap L^2(\Omega)\!\to\! \mathbb{R}}$, defined by\vspace*{-1.5mm}
	\begin{align}
		I(v):=\vert \textup{D}v\vert (\Omega)+\frac{\alpha}{2}\|v-g\|^2_{L^2(\Omega)}\label{ROF-primal}
	\end{align}
	for all $v\in BV(\Omega)\cap L^2(\Omega)$. In \cite[Theorem 10.5 \& Theorem 10.6]{Bar15}, it is established that for~every~${g\in L^2(\Omega)}$, there exists a unique minimizer $u\in BV(\Omega)\cap L^2(\Omega)$~of ${I\hspace*{-0.1em}:\hspace*{-0.1em}BV(\Omega)\cap L^2(\Omega)\hspace*{-0.1em}\to\hspace*{-0.1em} \mathbb{R}}$. If $g\hspace*{-0.1em}\in\hspace*{-0.1em} L^\infty(\Omega)$, then $u\hspace*{-0.1em}\in \hspace*{-0.1em}L^\infty(\Omega)$~with~${\|u\|_{L^\infty(\Omega)}\!\leq\! \|g\|_{L^\infty(\Omega)}}$ (cf. \cite[Proposition 10.2]{Bar15}). In \cite[Theorem 2.2]{HK04}, it is shown that the corresponding dual problem to \eqref{ROF-primal} determines a~vector~field~${z\!\in\! W^2_N(\textup{div};\Omega)\cap L^\infty(\Omega;\mathbb{R}^d)}$,~where $\Gamma_N\!=\!\partial\Omega$, that is maximal  for ${D\!:\!W^2_N(\textup{div};\Omega)\hspace*{-0.1em}\cap\hspace*{-0.1em} L^\infty(\Omega;\mathbb{R}^d)\!\to\! \mathbb{R}\hspace*{-0.1em}\cup\hspace*{-0.1em}\{-\infty\}}$,~defined~by\vspace*{-1mm}
	\begin{align}
		D(y):=-\frac{1}{2\alpha}\|\textup{div}(y)+g\|_{L^2(\Omega)}^2+\frac{\alpha}{2}\|g\|_{L^2(\Omega)}^2-I_{K_1(0)}(y)\label{ROF-dual}
	\end{align}
	for all $\smash{y\!\in\! W^2_N(\textup{div};\Omega)\hspace*{-0.1em}\cap\hspace*{-0.1em} L^\infty(\Omega;\mathbb{R}^d)}$, where $\smash{I_{K_1(0)}\!:\!L^\infty(\Omega;\mathbb{R}^d)\!\to\! \mathbb{R}\hspace*{-0.1em}\cup\hspace*{-0.1em}\{\infty\}}$~is~\mbox{defined} by $I_{K_1(0)}(y)\hspace*{-0.1em}:=\hspace*{-0.1em}0$ if $y\in L^\infty(\Omega;\mathbb{R}^d)$~with~${\|y\|_{L^\infty(\Omega;\mathbb{R}^d)}\leq 1}$~and~${I_{K_1(0)}(y):=\infty}$~else. Apart from that, in \cite[Theorem 2.2]{HK04}, it is shown that \eqref{ROF-dual} possesses a maximizer ${z\in W^2_N(\textup{div};\Omega)\cap L^\infty(\Omega;\mathbb{R}^d)}$, which satisfies the strong~duality~principle\vspace*{-1mm}
	\begin{align}
		I(u)=D(z).\label{strong_duality_principle}
	\end{align}\newpage
	\hspace*{-5mm}The \hspace*{-0.1mm}strong \hspace*{-0.1mm}duality \hspace*{-0.1mm}principle \hspace*{-0.1mm}\eqref{strong_duality_principle}, \hspace*{-0.1mm}appealing \hspace*{-0.1mm}to \hspace*{-0.1mm}\cite[\!Proposition \!10.4]{Bar15} \hspace*{-0.1mm}and~\hspace*{-0.1mm}\mbox{referring}~\hspace*{-0.1mm}to standard convex optimization arguments,~is~equivalent~to~the~optimality~relations
	\begin{align}
		\textup{div}(z)=\alpha (u-g)\quad\text{ in }L^2(\Omega),\qquad \vert \textup{D}u\vert(\Omega)=-(u,\textup{div}(z)).\label{optimality-relations}
	\end{align}

	\subsection{The discretized Rudin--Osher--Fatemi (ROF) model}\label{subsec:discrete_ROF_model}
	\qquad Given some $g\hspace*{-0.1em}\in \hspace*{-0.1em}L^2(\Omega)$ and $\alpha\hspace*{-0.1em}>\hspace*{-0.1em}0$, with $g_h\hspace*{-0.1em}:=\hspace*{-0.1em}\Pi_hg\hspace*{-0.1em}\in\hspace*{-0.1em} \mathcal{L}^0(\mathcal{T}_h)$, the discretized ROF \hspace*{-0.1mm}model \hspace*{-0.1mm}proposed \hspace*{-0.1mm}by \hspace*{-0.1mm}\cite{CP20} \hspace*{-0.1mm}determines \hspace*{-0.1mm}a \hspace*{-0.1mm}Crouzeix--Raviart~\hspace*{-0.1mm}\mbox{function}~\hspace*{-0.1mm}${u_h\!\in\! \mathcal{S}^{1,cr}(\mathcal{T}_h)}$ that~is~\mbox{minimal} for  $I_h:\mathcal{S}^{1,cr}(\mathcal{T}_h)\to \mathbb{R}$, defined by 
	\begin{align}
		I_h(v_h):=\|\nabla_hv_h\|_{L^1(\Omega;\mathbb{R}^d)}+\frac{\alpha}{2}\|\Pi_hv_h-g_h\|^2_{L^2(\Omega)}\label{discrete-ROF-primal}
	\end{align}
	for all $v_h\in \mathcal{S}^{1,cr}(\mathcal{T}_h)$. In \cite{CP20} and \cite{Bar21}, it has  been~shown~that~the~corresponding dual problem to \eqref{discrete-ROF-primal} determines a  Raviart--Thomas vector field ${z_h\in \mathcal{R}T^0_N(\mathcal{T}_h)}$, where ${\Gamma_N=\partial\Omega}$, that is maximal for  $D_h:\mathcal{R}T^0_N(\mathcal{T}_h)\to \mathbb{R}\cup\{-\infty\}$, defined by 
	\begin{align}
		D_h(y_h):=-\frac{1}{2\alpha}\|\textup{div}(y_h)+g_h\|_{L^2(\Omega)}^2+\frac{\alpha}{2}\|g_h\|_{L^2(\Omega)}^2-I_{K_1(0)}(\Pi_hy_h)\label{discrete-ROF-dual}
	\end{align}
	for all $y_h\!\in\! \mathcal{R}T^0_N(\mathcal{T}_h)$. Apart from that, in \cite{CP20} and \cite{Bar21}, it has been established~that a discrete weak~duality~principle holds, i.e., it holds
	\begin{align}
	\inf_{v_h\in \mathcal{S}^{1,cr}(\mathcal{T}_h)}{I_h(v_h)}\ge\sup_{y_h\in\mathcal{R}T^0_N(\mathcal{T}_h) }{ D_h(y_h)},\label{discrete_weak_duality_principle}
	\end{align}
	which is a cornerstone of the error analysis for \eqref{discrete-ROF-primal}. 
	In particular,~note~that~for~the validity of \eqref{discrete_weak_duality_principle} the $L^2$--projection operator~$\Pi_h$ in \eqref{discrete-ROF-primal} and \eqref{discrete-ROF-dual} plays~a~key~role.
	
	\subsection{Piece-wise Lipschitz, but not globally Lipschitz, continuous  solution to \eqref{ROF-dual}}\label{subsec:irregular_example}
	
	\qquad In \cite[Section 3]{BTW21}, the construction of an input data $g\in BV(\Omega)\cap L^\infty(\Omega)$ that leads to a solution $z\in W^\infty_N(\textup{div};\Omega)$ to \eqref{ROF-dual} such that ${z\notin W^{1,\infty}(\Omega;\mathbb{R}^2)}$,~in~essence, is based on the asymmetry of the function
	\begin{align}
		g:=\chi_{B_r^2(re_1)}-\chi_{B_r^2(-re_1)}\in BV(\Omega)\cap L^\infty(\Omega)\label{asym_data}
	\end{align}
	defined \hspace*{-0.1mm}on \hspace*{-0.1mm}a \hspace*{-0.1mm}domain \hspace*{-0.1mm}that \hspace*{-0.1mm}is \hspace*{-0.1mm}symmetric \hspace*{-0.1mm}with \hspace*{-0.1mm}respect \hspace*{-0.1mm}to \hspace*{-0.1mm}the \hspace*{-0.1mm}$\mathbb{R}e_2$--\hspace*{-0.1mm}axis.\!\footnote{For every $i=1,\dots,d$, we denote by $e_i\in \mathbb{S}^{d-1}$, the $i$--th. unit vector.}~\hspace*{-0.1em}More~\hspace*{-0.1mm}\mbox{precisely}, using this asymmetry property, it is possible to reduce the minimization problem \eqref{ROF-primal} on $\Omega$ into two independent minimization problems on $\Omega^+:=\Omega\cap (\mathbb{R}_{> 0}\times \mathbb{R})$ and $\Omega^-\hspace*{-0.1em}:=\hspace*{-0.1em}\Omega\cap (\mathbb{R}_{< 0}\hspace*{-0.1em}\times\hspace*{-0.1em} \mathbb{R})$ for which explicit solutions ${u^{\pm}\hspace*{-0.1em}\in\hspace*{-0.1em} BV(\Omega^{\pm})\cap L^\infty(\Omega^{\pm})}$~exist. In this way, the following result could be derived.
	
	\begin{proposition}\label{asym_primal_solution}
		Let $\Omega\subseteq \mathbb{R}^2$ be symmetric with respect to the $\mathbb{R}e_2$--axis and let $r>0$ be such that ${B_r^2(\pm r e_1)\subset\subset \Omega}$. Then, for \eqref{asym_data}~and~${\alpha>0}$, the minimizer~of ${I\hspace*{-0.05em}:\hspace*{-0.05em}\mathcal{A}\hspace*{-0.05em}\to\hspace*{-0.05em} \mathbb{R}}$, where ${\mathcal{A}\hspace*{-0.05em}:=\hspace*{-0.05em}\{u\hspace*{-0.05em}\in\hspace*{-0.05em} BV(\Omega)\hspace*{-0.05em}\cap\hspace*{-0.05em} L^\infty(\Omega)\mid\textup{tr}(u)\hspace*{-0.05em}=\hspace*{-0.05em}0\textup{ in }L^1(\partial\Omega) \}}$,~is~given~via 
		\begin{align}
			u=\max\Big\{0,1-\frac{2}{\alpha r}\Big\} g\in \mathcal{A}.\label{eq:asym_primal_solution}
		\end{align}
	\end{proposition}
	
	\begin{proof}
		See \cite[Proposition~3.1]{BTW21}.
	\end{proof}

	\hspace*{-1mm}Combining \hspace*{-0.1mm}the \hspace*{-0.1mm}representation \hspace*{-0.1mm}formula \hspace*{-0.1mm}\eqref{eq:asym_primal_solution} \hspace*{-0.1mm}and \hspace*{-0.1mm}the \hspace*{-0.1mm}optimality~\hspace*{-0.1mm}conditions~\hspace*{-0.1mm}\eqref{optimality-relations}, it turns out that there exists no Lipschitz continuous dual solution to the setting described in Proposition \ref{asym_primal_solution}.
	
	\begin{corollary}\label{asym_primal_solution_corollary}
		Let the assumptions of Proposition \ref{asym_primal_solution} be satisfied~with~${\alpha r>2}$. Then,  any dual solution $z\in W^\infty(\textup{div};\Omega)$ to \eqref{eq:asym_primal_solution} is not $\theta$--Hölder~continuous if ${\theta> \frac{1}{2}}$.
	\end{corollary}

	\begin{proof}
	See \cite[Corollary 3.2]{BTW21}.
	\end{proof}

An example of  a not Lipschitz continuous dual solution to \eqref{eq:asym_primal_solution}~is~the~\mbox{following},
which is separately Lipschitz continuous on $\Omega^+$ and $\Omega^-$, resp., and jumps over the $\mathbb{R}e_2$--axis. We will resort to this dual solution to derive optimal error estimates for the setting described in Proposition \ref{asym_primal_solution}, which has~already~been~reported~in \cite[Example~6.1]{BTW21}.
	
	\begin{proposition}\label{asym_dual_solution}
		Let $\Omega\subseteq \mathbb{R}^2$  and $r>0$ be  such as in Proposition \ref{asym_primal_solution}. Moreover, let $\alpha>0$ be such that $\alpha r>2$. Then, the vector field $z:\Omega\subseteq \mathbb{R}^2\to \mathbb{R}^2$, defined by
		\begin{align*}
			z(x):=\begin{cases}
				\mp\tfrac{1}{r}(x\mp re_1)&\quad\text{ if }\vert x\mp re_1\vert <r\\
				\mp\tfrac{r}{\vert x\mp re_1\vert^2}(x\mp re_1)&\quad\text{ if }\vert x\mp re_1\vert\ge r
			\end{cases}
		\end{align*}\vspace*{1mm}
		for all $x\!\in\! \Omega$, satisfies \hspace*{-0.075em}$z\!\in\! W^\infty(\textup{div};\hspace*{-0.1em}\Omega)$, \hspace*{-0.075em}${\|z\|_{L^\infty(\Omega;\mathbb{R}^d)}\!\leq\! 1}$, \hspace*{-0.075em}$\vert \textup{D}u\vert(\Omega)\!=\!-(u,\textup{div}(z))_{L^2(\Omega)}$ and \hspace*{-0.075em}$ \textup{div}(z)\!=\!\alpha(u-g)$ \hspace*{-0.075em}in \hspace*{-0.075em}$L^\infty(\Omega)$, \hspace*{-0.075em}where~\hspace*{-0.075em}${u\!\in \!\mathcal{A}}$~\hspace*{-0.075em}is~\hspace*{-0.075em}defined~\hspace*{-0.075em}by\hspace*{-0.075em}~\eqref{eq:asym_primal_solution},~\hspace*{-0.075em}i.e.,~\hspace*{-0.075em}${z\!\in\! W^\infty(\textup{div};\Omega)}$ is a dual solution to \eqref{eq:asym_primal_solution}.
	\end{proposition}
	
	\begin{proof}
		Apparently, we have that $z\hspace*{-0.1em}\in\hspace*{-0.1em}  L^\infty(\Omega;\mathbb{R}^d)$ with ${\|z\|_{L^\infty(\Omega;\mathbb{R}^d)}\hspace*{-0.1em}\leq \hspace*{-0.1em}1}$.~In~\mbox{addition}, 
		it is not difficult to see~that  $\smash{z|_{\Omega^{\pm}}\in W^{1,\infty}(\Omega^{\pm};\mathbb{R}^d)}$. Since ${z|_{\Omega^+}\cdot n_{\Omega^+}\hspace*{-0.1em}=\hspace*{-0.1em}-z|_{\Omega^-}\cdot n_{\Omega^-}}$ on $\mathbb{R}e_2\cap\Omega$, we find~that $z\!\in\! W^2(\textup{div};\Omega)$. \!It is well-known,~cf.~\mbox{\cite[Example~\!10.4]{Bar15}},~that
		\begin{align*}
			\vert \textup{D}u\vert(\Omega^{\pm})=-(u,\textup{div}(z))_{L^2(\Omega^{\pm})},\qquad \textup{div}(z)=\alpha(u-g)\quad\text{ in } L^\infty(\Omega^{\pm}).
		\end{align*}
		Thus, we have that $ \textup{div}(z)=\alpha(u-g)$ in $L^\infty(\Omega)$, which implies that ${z\!\in\! W^\infty(\textup{div};\Omega)}$.
		 Apart from that, using that $u=0$ continuously in $\mathbb{R}e_2\cap\Omega$, we finally~conclude that $\vert \textup{D}u\vert(\Omega)=-(u,\textup{div}(z))_{L^2(\Omega)}$, i.e., $z\in W^\infty(\textup{div};\Omega)$ is a dual solution~to~\eqref{eq:asym_primal_solution}.\vspace*{5mm}
	\end{proof}

	\section{Error estimates depending on Sobolev regularity}\label{sec:lipschitz_truncation}
	
\qquad 	The validity of quasi-optimal error estimates for the finite element approxima-tion of total-variation regularized minimization problems~by~means~of~the~Crouzeix--Raviart element in the case of an existing Lipschitz continuous solution~to~\eqref{ROF-dual}~in \cite{CP20,Bar21}, \hspace*{-0.15mm}in \hspace*{-0.15mm}essence, \hspace*{-0.15mm}is \hspace*{-0.15mm}based \hspace*{-0.15mm}on \hspace*{-0.15mm}four \hspace*{-0.15mm}results: \!The \hspace*{-0.15mm}discrete \hspace*{-0.15mm}weak~\hspace*{-0.15mm}duality~\hspace*{-0.15mm}\mbox{principle~\!\eqref{discrete_weak_duality_principle}}, the discrete strong coercivity~of~$I_h:\mathcal{S}^{1,cr}(\mathcal{T}_h)\to \mathbb{R}$,~i.e.,
	\begin{align}
		\frac{\alpha}{2}\| \Pi_h(v_h-u_h)\|_{L^2(\Omega)}^2\leq I_h(v_h)-I_h(u_h)\label{discrete_strong_coercivity}
	\end{align}\vspace*{1mm}
	for all $\smash{v_h\in \mathcal{S}^{1,cr}(\mathcal{T}_h)}$, where $\smash{u_h\in \mathcal{S}^{1,cr}(\mathcal{T}_h)}$ is the minimum of  $\smash{I_h:\mathcal{S}^{1,cr}(\mathcal{T}_h)\to \mathbb{R}}$, the strong duality principle  \eqref{strong_duality_principle}, and the existence of appropriate primal~and~dual quasi-interpolants, guaranteed through the following two lemmas:
	
For the benefit of readability and without loss of generality, we~assume~for~the remainder of this article, if not otherwise specified, that $\alpha=1$.\newpage
	\begin{lemma}[Primal quasi-interpolant]\label{primal_quasi_interpolant0}
		For every $u\in BV(\Omega)\cap L^\infty(\Omega)$, there exists  a Crouzeix--Raviart function $\tilde{u}_h\in \mathcal{S}^{1,cr}(\mathcal{T}_h)$ with the following properties:
		\begin{description}[noitemsep,topsep=2pt,labelwidth=\widthof{P.3},leftmargin=!,font=\normalfont\itshape]
			\item[(P.1)]\hypertarget{(P.1)}{}
			$\|\nabla_h \tilde{u}_h\|_{L^1(\Omega;\mathbb{R}^d)}\leq \vert \textup{D}u\vert(\Omega)$,
			\item[(P.2)]\hypertarget{(P.2)}{}
			$\|u-\tilde{u}_h\|_{L^1(\Omega)}\leq c_{cr}h\vert \textup{D}u\vert(\Omega)$,
			\item[(P.3)] \hypertarget{(P.3)}{}
			$\|\tilde{u}_h\|_{L^\infty(\Omega)}\leq c_d\|u\|_{L^\infty(\Omega)}$,
			\item[(P.4)] \hypertarget{(P.4)}{}
			$I_h(\tilde{u}_h)\leq I(u)+2c_dc_{cr}\|u\|_{L^\infty(\Omega)}\vert \textup{D}u\vert(\Omega)h-\frac{1}{2}\smash{\|g-g_h\|_{L^2(\Omega)}^2}$.
		\end{description}
	\end{lemma}
	
		\begin{proof}
		See \cite[Lemma 4.4]{Bar21} or \cite[Section 5]{CP20}.
	\end{proof}

	\begin{lemma}[Dual quasi--interpolant] \label{dual_quasi_interpolant0} For every ${z\in W^{1,\infty}(\Omega;\mathbb{R}^d)\cap W^2_N(\textup{div};\Omega)}$ such that $\| z\|_{L^\infty(\Omega;\mathbb{R}^d)}\!\leq\! 1$, there exists a Raviart--Thomas vector field ${\tilde{ z}_h\!\in\! \mathcal{R}T^0_N(\mathcal{T}_h)}$ with the following properties:
		\begin{description}[noitemsep,topsep=2pt,labelwidth=\widthof{(D)},leftmargin=!,font=\normalfont\itshape]
			\item[(D.1)] \hypertarget{(D.1)}{} $\|\Pi_h\tilde{ z}_h\|_{L^\infty(\Omega;\mathbb{R}^d)}\leq 1$,
			
			\item[(D.2)] \hypertarget{(D.2)}{} $D_h(\tilde{ z}_h)\!\ge\! D( z)-c_{\mathcal{R}T}\|\nabla z\|_{L^\infty(\Omega;\mathbb{R}^{d\times d})}\|g\|_{L^2(\Omega)}\|\textup{div}(z)\|_{L^2(\Omega)}h-\frac{1}{2}\|g-g_h\|_{L^2(\Omega)}^2$.
		\end{description}
	\end{lemma}

	\begin{proof}
		See \cite[Lemma 4.5]{Bar21} or \cite[Section 5]{CP20}.
	\end{proof}
	
	While Lemma \ref{primal_quasi_interpolant0} does not impose restrictive assumptions on the minimum $u\in BV(\Omega)\cap L^2(\Omega)$, since already $u\in L^\infty(\Omega)$ if $g\in L^\infty(\Omega)$ (cf. \cite[Proposition~10.2]{Bar15}), the required Lipschitz continuity of a solution $z\!\in\! W_N^\infty\hspace*{-0.1em}(\textup{div};\hspace*{-0.1em}\Omega\hspace*{-0.05em})$~to~\hspace*{-0.1em}\eqref{ROF-dual}~in~Lemma~\hspace*{-0.1em}\ref{dual_quasi_interpolant0} is often not fulfilled, cf. \!\cite{BTW21} or Section \ref{subsec:irregular_example}.
	 We resort to the Lipschitz~\mbox{truncation} technique to fulfill the Lipschitz continuity requirement~on~a~solution~to~\eqref{ROF-dual}~in Lemma \ref{dual_quasi_interpolant0} at least in an approximate sense and, in this way, derive error estimates that depend directly on the Sobolev regularity of a solution~$\hspace*{-0.1em}{z\!\in\! W^{1,p}(\Omega;\hspace*{-0.1em}\mathbb{R}^d)}\hspace*{-0.1em}$~to~\hspace*{-0.1em}\eqref{ROF-dual}. 
	 The main advantage of this approach is that the Lipschitz truncation technique is based on local arguments, while regularization~by~\hspace*{-0.1em}\mbox{convolution}~as~in~\mbox{\cite[\!Section~\!5.2]{CP20}}, for example, operates highly non-local and, therefore, wipes out point-wise and/or local properties of a solution to \eqref{ROF-dual} that potentially could have~been~\mbox{incorporated}. 
	 To be more precise, in \cite[Section~5.2]{CP20},  the requirement  $g\in BV(\Omega)\cap L^\infty(\Omega)$~was needed  to estimate  $\|\textup{div}(z)-\textup{div}(z_\varepsilon)\|_{\smash{L^1(\Omega)}}$, where $z_\varepsilon\hspace*{-0.1em} :=\hspace*{-0.1em} z\circ \omega_\varepsilon\hspace*{-0.1em}\in\hspace*{-0.1em} C^\infty(\mathbb{R}^d;\mathbb{R}^d)$,~${\varepsilon\hspace*{-0.1em}>\hspace*{-0.1em}0}$, denotes the convolution with a suitably scaled kernel $\omega_\varepsilon\in C^\infty_0(\mathbb{R}^d)$,~by~$\varepsilon\vert \textup{D}g\vert(\Omega)$. In contrast to that, if $z_\lambda \in W^{1,\infty}(\mathbb{R}^d;\mathbb{R}^d)$, $\lambda>0$, denotes the Lipschitz truncation of a suitable extension $\overline{z}\in W^{1,p}(\mathbb{R}^d;\mathbb{R}^d)$ of ${z\in W^{1,p}(\Omega;\mathbb{R}^d)}$, then we can exploit the particular properties	$\nabla z_\lambda\!=\!\nabla \overline{z}$ in $\{z_\lambda = \overline{z}\}$ and ${\vert \{z_\lambda \neq \overline{z}\}\vert \!\leq\! \vert \{\mathcal{M}(\nabla \overline{z})>\lambda \}\vert}$\footnote{Here, $\mathcal{M}\!:\!L^p(\mathbb{R}^d;\mathbb{R}^l)\!\to\! L^p(\mathbb{R}^d;\mathbb{R}^l)$, $d,l\in \mathbb{N}$, defined by $\smash{\mathcal{M}(f)(x)\!:=\!\sup_{r>0}{\fint_{B_r^d(x)}{\vert f(y)\vert \,\textup{d}y}}}$ for a.e. $x\in \mathbb{R}^d$ and all $f\in L^p(\mathbb{R}^d;\mathbb{R}^l)$,~denotes~the~\mbox{Hardy--Littlewood--Maximal operator}.}, to~conclude~that $\|\textup{div}(z)-\textup{div}(z_\lambda)\|_{\smash{L^1(\Omega)}}
	 \leq c\lambda ^{1-p}\|\nabla z\|_{L^p(\Omega;\mathbb{R}^{d\times d})}$. 
	In~this~way,~we obtain the same  rate $\smash{\mathcal{O}(h^{\frac{1}{4}})}$ in \cite[Section~5.2]{CP20} without 
	the assumption $g\in BV(\Omega)$ but need to require $\smash{z\in W^{1,3}(\Omega;\mathbb{R}^d)}\cap L^\infty(\Omega;\mathbb{R}^d)$~instead~of only
	$\smash{z\in W^\infty_N(\textup{div};\Omega)}$.

	\begin{theorem}[Lipschitz truncation technique]\label{lipschitz_truncation_technique}
		Let $z\!\in\! W^{1,p}(\mathbb{R}^d;\mathbb{R}^d)$, ${p\!\in\! \left[1,\infty\right)}$, and $\theta,\lambda\!>\!0$. Then, there is a Lipschitz continuous~vector~field~${z_{\theta,\lambda}\!\in\! W^{1,\infty}(\mathbb{R}^d;\mathbb{R}^d)}$ and a constant ${c_{\textup{LT}}\hspace*{-0.1em}>\hspace*{-0.1em}0}$, which does not depend on ${p\hspace*{-0.1em}\in\hspace*{-0.1em} \left[1,\infty\right)}$~and~${\theta,\lambda\hspace*{-0.1em}>\hspace*{-0.1em}0}$,~such~that the following statements apply:
		\begin{description}[noitemsep,topsep=2pt,labelwidth=\widthof{(LT.4)},leftmargin=!,font=\normalfont\itshape]
			\item[(LT.1)]\hypertarget{LT.1}{}
			$\|z_{\theta,\lambda}\|_{L^\infty(\mathbb{R}^d;\mathbb{R}^d)}\leq \theta$,
			\item[(LT.2)]\hypertarget{LT.2}{}
			$\|\nabla z_{\theta,\lambda}\|_{L^\infty(\mathbb{R}^d;\mathbb{R}^{d\times d})}\leq c_{\textup{LT}}\lambda$,
			\item[(LT.3)]\hypertarget{LT.3}{}
			 $\vert\{z_{\theta,\lambda}\neq z\}\vert\leq \vert\{\mathcal{M}(z)>\theta\}\vert+\vert\{ \mathcal{M}(\nabla z)> \lambda\}\vert$,
			\item[(LT.4)] \hypertarget{(LT.4)}{}
			$\nabla z_{\theta,\lambda}=\nabla z$ in $\{z_{\theta,\lambda}= z\}$. 
		\end{description}
		
	\end{theorem}
	
	\begin{proof}
		See the first part of the proof of \cite[Theorem 2.3]{DMS08} or \cite[Section 1.3.3]{MZ97}.
	\end{proof}

	A crucial property of the Lipschitz truncation technique for this article~is~that, similar to regularization by convolution, it does not increase the maximal~length~of a vector  field.
	
	\begin{remark}[Maximal \hspace*{-0.15mm}length \hspace*{-0.15mm}preservation \hspace*{-0.15mm}of \hspace*{-0.15mm}the \hspace*{-0.15mm}Lipschitz \hspace*{-0.15mm}truncation \hspace*{-0.15mm}technique]\label{modulus_preservation}
		If ${z\in W^{1,p}(\mathbb{R}^d;\mathbb{R}^d)\cap L^\infty(\mathbb{R}^d;\mathbb{R}^d)}$, ${p\in \left[1,\infty\right)}$, for some $\theta>0$ has~the~property $\|z\|_{L^\infty(\mathbb{R}^d;\mathbb{R}^d)} \leq \theta$, then
		\begin{align*}
			\mathcal{M}(z)(x)=\sup_{r>0}{\fint_{B_r^d(x)}{\vert z(y)\vert\,dy}}\leq \| z\|_{L^\infty(\mathbb{R}^d;\mathbb{R}^d)} \leq \theta\quad\text{ for a.e. }x\in \mathbb{R}^d,
		\end{align*}
		i.e., $\vert\{ \mathcal{M}(  z)> \theta\}\vert=0$, which (cf. Theorem~\ref{lipschitz_truncation_technique},~(\hyperlink{(LT.3)}{LT.3}))  for arbitrary~${\lambda>0}$~yields 
		\begin{align}
			\vert\{ z_{\theta,\lambda}\neq  z\}\vert\leq \vert\{ \mathcal{M}(\nabla  z)> \lambda\}\vert.\label{rmk1.1}
		\end{align}
		Through the combination of the $p$--type Tschebyscheff--Markoff--inequality, i.e.,
		$\smash{\vert\{\mathcal{M}(\nabla z)\!>\!\lambda\}\vert\!\leq\! \lambda^{-p}\|\mathcal{M}(\nabla z)\|_{L^{p}(\mathbb{R}^d;\mathbb{R}^{d\times d})}^p}$,  and the strong type $(p,p)$--estimate of the Hardy--Littlewood--Maximal operator (cf. \!\cite[\!Theorem \!1.22]{MZ97}),~i.e.,~for~${c_{\mathcal{M}}\!>\!0}$\footnote{More precisely, one has $c_{\mathcal{M}}\!=\!2 \big(\frac{p}{p-1}\big)^{\frac{1}{p}}5^{\frac{d}{p}}$, implying the limit behavior $c_{\mathcal{M}}\!\to\! 2$ for $(p \!\to\! \infty)$\vspace*{-3mm}.}\!, $\smash{\|\mathcal{M}(\nabla z)\|_{L^{p}(\mathbb{R}^d;\mathbb{R}^{d\times d})}^p\leq c _{\mathcal{M}} \|\nabla z\|_{L^{p}(\mathbb{R}^d;\mathbb{R}^{d\times d})}^p}$,
		we deduce from \eqref{rmk1.1} that
		\begin{align}
			\vert\{ z_{\theta,\lambda}\neq  z\}\vert\leq c_{\mathcal{M}} \lambda^{-p}\|\nabla z\|_{L^{p}(\mathbb{R}^d;\mathbb{R}^{d\times d})}^p.\label{rmk1.2}
		\end{align}
		By means of \eqref{rmk1.2}, also using Theorem~\ref{lipschitz_truncation_technique},~(\hyperlink{(LT.2)}{LT.2}) \& (\hyperlink{(LT.4)}{LT.4}), we, then,~deduce~that
		\begin{align}
				\|\nabla z_{\theta,\lambda}\|_{L^{p}(\mathbb{R}^d;\mathbb{R}^{d\times d})}&=	\|\nabla z\chi_{\{ z_{\theta,\lambda}=  z\}}\|_{L^{p}(\mathbb{R}^d;\mathbb{R}^{d\times d})}+\|\nabla z_{\theta,\lambda}\chi_{\{ z_{\theta,\lambda}\neq  z\}}\|_{L^{p}(\mathbb{R}^d;\mathbb{R}^{d\times d})}\notag
				\\&\leq\|\nabla z\|_{L^{p}(\mathbb{R}^d;\mathbb{R}^{d\times d})}+c_{\textup{LT}}\lambda\vert \{ z_{\theta,\lambda}\neq  z\}\vert^{\smash{\frac{1}{p}}}\label{rmk1.3}
				\\&\leq \big(1+{c_{\mathcal{M}}}^{\smash{\!\frac{1}{p}}}c_{\textup{LT}}\big) \|\nabla z\|_{L^{p}(\mathbb{R}^d;\mathbb{R}^{d\times d})}.\notag
		\end{align}
	\end{remark}

	Through the combination of Lemma \ref{dual_quasi_interpolant0}, Theorem \ref{lipschitz_truncation_technique} and Remark \ref{modulus_preservation}, we arrive at the following result providing an admissible dual quasi-interpolant whose particular \hspace*{-0.1mm}properties \hspace*{-0.1mm}depend \hspace*{-0.1mm}directly \hspace*{-0.1mm}on \hspace*{-0.1mm}the \hspace*{-0.1mm}Sobolev~\hspace*{-0.1mm}regularity~\hspace*{-0.1mm}of~\hspace*{-0.1mm}a~\hspace*{-0.1mm}\mbox{solution}~\hspace*{-0.1mm}to~\hspace*{-0.1mm}\eqref{ROF-dual}.
	
	\begin{lemma}[Dual quasi--interpolant depending on Sobolev regularity~for~${\Gamma_N\hspace*{-0.15em}=\hspace*{-0.15em}\emptyset}$]\label{dual_quasi_interpolant1}
		Let $g\in L^\infty(\Omega)$ and let  ${z\in W^{1,p}(\Omega;\mathbb{R}^d)\cap W^\infty(\textup{div};\Omega)}$, ${p\in\left[2,\infty\right)}$, be such that  $\|  z\|_{L^\infty(\Omega;\mathbb{R}^d)}\leq 1$.  Then, there exists a Raviart--Thomas vector field  ${\tilde{ z}_h\in \mathcal{R}T^0(\mathcal{T}_h)}$ with the following properties:
		\begin{description}[noitemsep,topsep=2pt,labelwidth=\widthof{(D.2)},leftmargin=!,font=\normalfont\itshape]
				\item[($D_p$.1)] \hypertarget{(Dp.1)}{} $\|\Pi_h\tilde{ z}_h\|_{L^\infty(\Omega;\mathbb{R}^d)}\leq 1$.
				
				\item[($D_p$.2)] \hypertarget{(Dp.2)}{} $\smash{D_h(\tilde{ z}_h)\ge D( z)-c_ph^{\frac{p-2}{p-1}}-\frac{1}{2}\|g-g_h\|_{L^2(\Omega)}^2}$, where
				\begin{align*}c_p(z)&:=2c_{\mathcal{R}T}\|g\|_{L^2(\Omega)}d^{\frac{1}{2}}\big(1+{c_{\mathcal{M}}}^{\smash{\!\frac{1}{2}}}c_{\textup{LT}}\big) c_{\textup{E}}\|\nabla z\|_{L^2(\Omega;\mathbb{R}^{d\times d})}c_{\textup{LT}}\\&\quad+8dc_{\textup{LT}}^2c_{\mathcal{M}}c_{\textup{E}}^p\|\nabla z\|_{L^p(\Omega;\mathbb{R}^{d\times d})}^p.
				\end{align*}
				Here, $c_{\textup{E}}>0$ is the Lipschitz constant of the lower-order extension~operator $P\hspace*{-0.05em}:\hspace*{-0.05em}W^{1,q}(\Omega;\mathbb{R}^l)\hspace*{-0.05em}\to \hspace*{-0.05em}W^{1,q}(\mathbb{R}^d;\mathbb{R}^l)$, $q\hspace*{-0.05em}\in \hspace*{-0.05em}[1,\infty]$, constructed~in~\mbox{\cite[\!Section~\!9.2]{Bre10}}, which does not depend on $q \in [1,\infty]$.
		\end{description}
	\end{lemma}

	\begin{remark}
		\begin{description}[noitemsep,topsep=2pt,labelwidth=\widthof{(iii)},leftmargin=!,font=\normalfont\itshape]
			\item[(i)]
			
			\hspace*{-2mm}The arguments remain valid for $\Gamma_D\!\neq \! \partial \Omega$ if for  $z\!\in\! W^{1,p}(\Omega;\mathbb{R}^d)$ $\cap W^\infty_N(\textup{div};\Omega)$, $p\!\in\!\left[2,\infty\right)$, such that $\|  z\|_{L^\infty(\Omega;\mathbb{R}^d)}\!\leq\! 1$,~there~exists~an~\mbox{extension}  $\smash{\overline{ z}\in W^{1,p}(\mathbb{R}^d;\mathbb{R}^d)}$ with $\overline{ z}|_{\Omega}\!=\! z$~and~$\smash{\|\overline{ z}\|_{L^\infty(\mathbb{R}^d;\mathbb{R}^d)}\!\leq\! 1}$ and~if~for~this~\mbox{extension}, the Lipschitz truncation $\overline{z}_{1,\lambda}\in W^{1,\infty}(\mathbb{R}^d;\mathbb{R}^d)$ from Theorem \ref{lipschitz_truncation_technique} satisfies $ \overline{z}_{1,\lambda}\cdot n =0$ in $\Gamma_N$.
				
			\item[(ii)] In general, the constant $c_p>0$ deteriorates as $p\to \infty$,~i.e.,~${c_p\to \infty}$~${(p\to\infty)}$.
		\end{description}
	\end{remark}
	
	\begin{proof} (of Lemma \ref{dual_quasi_interpolant1})
		Resorting to a lower-order extension operator, as, e.g., in \cite[Theorem~9.7]{Bre10}, we get some $\smash{\overline{ z}\!\in\! W^{1,p}(\mathbb{R}^d;\mathbb{R}^d)}$ with $\overline{ z}|_{\Omega}\!=\! z$~and~$\smash{\|\overline{ z}\|_{L^\infty(\mathbb{R}^d;\mathbb{R}^d)}\!\leq\! 1}$. Denote by $ z_{\lambda}\hspace*{-0.1em}:=\hspace*{-0.1em}\overline{ z}_{1,\lambda}\hspace*{-0.1em}\in\hspace*{-0.1em} W^{1,\infty}(\mathbb{R}^d;\mathbb{R}^d)$,~i.e., for $\theta\hspace*{-0.1em}=\hspace*{-0.1em}1$, the Lipschitz truncation~of   $\overline{z}\hspace*{-0.15em}\in\hspace*{-0.2em} W^{1,p}(\mathbb{R}^d;\mathbb{R}^d)$  in the sense of Theorem~\ref{lipschitz_truncation_technique}. \hspace*{-0.25em}Then, also using Remark~\ref{modulus_preservation},~we~get:
		\begin{description}[noitemsep,topsep=1pt,labelwidth=\widthof{(iii)},leftmargin=!,font=\normalfont\itshape]
			\item[($\alpha$)]
			$\| z_{\lambda}\|_{L^\infty(\mathbb{R}^d;\mathbb{R}^d)}\leq 1$,
			\item[($\beta$)]
			$\|\nabla z_{\lambda}\|_{L^\infty(\mathbb{R}^d;\mathbb{R}^{d\times d})}\leq c_{\textup{LT}}\lambda$,
			\item[($\gamma$)] $\vert\{ z_{\lambda}\neq  \overline{z}\}\vert\leq\vert\{ \mathcal{M}(\nabla  \overline{z})> \lambda\}\vert$,
			\item[($\delta$)] $\nabla z_{\lambda}=\nabla \overline{z}$ in $\{ z_{\lambda}=  \overline{z}\}$. 
		\end{description}
		For $ z_{\lambda}|_{\Omega}\in W^{1,\infty}(\Omega;\mathbb{R}^d)$ we obtain, in analogy with Lemma \ref{dual_quasi_interpolant0}, i.e., introducing
		$\smash{z_h^{\lambda}\hspace*{-0.05em}:=\hspace*{-0.05em}(\gamma_h^\lambda)^{-1}I_{\mathcal{R}T}z_{\lambda}\hspace*{-0.05em}\in\hspace*{-0.05em} \mathcal{R}T^0(\mathcal{T}_h)}$, where we define ${\gamma_h^\lambda\hspace*{-0.05em}:=\hspace*{-0.05em}1+c_{\mathcal{R}T}\|\nabla z_\lambda\|_{L^\infty(\Omega;\mathbb{R}^{d\times d})}h}$,
		 a dual quasi-interpolant  $\smash{z_h^{\lambda}\in \mathcal{R}T^0(\mathcal{T}_h)}$ such that both $\smash{\|  \Pi_h z_h^{\lambda}\|_{L^\infty(\Omega;\mathbb{R}^d)} \leq 1}$ and
		\begin{align}
			\begin{aligned}
				D_h( z_h^\lambda)&\ge D( z_\lambda)-c_{\mathcal{R}T}\|g\|_{L^2(\Omega)}\|\textup{div}( z_\lambda)\|_{L^2(\Omega)}\|\nabla z_\lambda\|_{L^\infty(\Omega;\mathbb{R}^{d\times d})}h\\&\quad-\tfrac{1}{2}\|g-g_h\|_{L^2(\Omega)}^2.	\end{aligned}\label{lem1.1}
		\end{align}
		Then, on the basis of \textit{($\beta$)} and \eqref{rmk1.3}, we find that
		\begin{align}
			\hspace*{-0.3cm}\begin{aligned}
				\|\textup{div}( z_\lambda)\|_{L^2(\Omega)}\|\nabla z_\lambda\|_{L^\infty(\Omega;\mathbb{R}^{d\times d})}&\leq 
				d^{\frac{1}{2}}\|\nabla z_\lambda\|_{L^2(\mathbb{R}^d;\mathbb{R}^{d\times d})}\|\nabla z_\lambda\|_{L^\infty(\mathbb{R}^d;\mathbb{R}^{d\times d})}
				\\&\leq 
				d^{\frac{1}{2}}\big(1+{c_{\mathcal{M}}}^{\smash{\!\frac{1}{2}}}c_{\textup{LT}}\big) \|\nabla \overline{z}\|_{L^2(\mathbb{R}^d;\mathbb{R}^{d\times d})}c_{\textup{LT}}\lambda
				\\&\leq d^{\frac{1}{2}}\big(1+{c_{\mathcal{M}}}^{\smash{\!\frac{1}{2}}}c_{\textup{LT}}\big) c_{\textup{E}} \|\nabla z\|_{L^2(\Omega;\mathbb{R}^{d\times d})}c_{\textup{LT}}\lambda.\end{aligned}\label{lem1.2}
		\end{align} 
		Using \textit{($\beta$)},  \textit{($\delta$)} and \eqref{rmk1.2}, also assuming that ${d^{\frac{1}{2}}c_{\textup{LT}}\lambda\hspace*{-0.15em}>\hspace*{-0.15em}\|\textup{div}( z)\|_{L^\infty(\Omega)}\hspace*{-0.15em}+\hspace*{-0.15em}2\|g\|_{L^\infty(\Omega)}}$, we further deduce that
		\begin{align}
				\vert D( z)-&D( z_{\lambda})\vert\leq \|(\textup{div}( z_{\lambda})\hspace*{-0.1em}-\hspace*{-0.1em}\textup{div}( z))\chi_{\{ z_\lambda\neq \overline{z} \}\cap \Omega}\|_{L^1(\Omega)}\|\textup{div}( z_{\lambda})\hspace*{-0.1em}+\hspace*{-0.1em}\textup{div}( z)\hspace*{-0.1em}-\hspace*{-0.1em}2g\|_{L^\infty(\Omega)}\notag
			\\\notag&\leq(d^{\frac{1}{2}}c_{\textup{LT}}\lambda\hspace*{-0.1em}+\hspace*{-0.1em}\|\textup{div}( z)\|_{L^\infty(\Omega)}\hspace*{-0.1em}+\hspace*{-0.1em}2\|g\|_{L^\infty(\Omega)})^2\vert\{ z_\lambda\neq \overline{z} \}\vert
				\\\notag&\leq\smash{4dc_{\textup{LT}}^2}\lambda^2 c_{\mathcal{M}}\lambda^{-p}\|\nabla \overline{z}\|_{L^p(\mathbb{R}^d;\mathbb{R}^{d\times d})}^p\notag
				\\&\leq \smash{4dc_{\textup{LT}}^2c_{\mathcal{M}}c_{\textup{E}}^p\|\nabla z\|_{L^p(\Omega;\mathbb{R}^{d\times d})}^p}\lambda^{2-p}.
	\label{lem1.3}
		\end{align}
		Therefore, on combining \eqref{lem1.1}, \eqref{lem1.2} and \eqref{lem1.3},  we observe that
		\begin{align}
			D_h( z_h^\lambda)&\ge D( z)-4dc_{\textup{LT}}^2c_{\mathcal{M}}c_{\textup{E}}^p\|\nabla z\|_{L^p(\Omega;\mathbb{R}^{d\times d})}^p\lambda^{2-p}-\tfrac{1}{2}\|g-g_h\|_{L^2(\Omega)}^2\label{lem1.4}\\&\quad-\smash{c_{\mathcal{R}T}\|g\|_{L^2(\Omega)}d^{\frac{1}{2}}\big(1+{c_{\mathcal{M}}}^{\frac{1}{2}}c_{\textup{LT}}\big) c_{\textup{E}}\|\nabla z\|_{L^2(\Omega;\mathbb{R}^{d\times d})}c_{\textup{LT}}\lambda h}.\notag
		\end{align}
		For $c_p>0$ defined as above and
		$\lambda = h^{-s}$, where $s>0$ is arbitrary,~\eqref{lem1.4}~yields
		\begin{align}
			D_h( z_h^\lambda)\ge D( z)-\tfrac{c_p}{2}(h^{s(p-2)}+h^{1-s})-\tfrac{1}{2}\|g-g_h\|_{L^2(\Omega)}^2.\label{lem1.5}
		\end{align}
		We have that $s(p-2)\hspace*{-0.1em}=\hspace*{-0.1em}1-s$ if and only if $\smash{s=\frac{1}{p-1}\hspace*{-0.1em}=\hspace*{-0.1em}\frac{p'}{p}}$. 
		Thus, for ${\lambda\hspace*{-0.1em} =\hspace*{-0.1em} h^{-s}}$,~$\smash{s\hspace*{-0.1em}=\hspace*{-0.1em}\frac{1}{p-1}}$ and $\smash{\tilde{z}_h:=z_h^\lambda\in \mathcal{R}T^0(\mathcal{T}_h)}$, 
		from \eqref{lem1.5}, it follows that~both~\textit{(\hyperlink{(Dp.1)}{$D_p$.1})}~and~\textit{(\hyperlink{(Dp.2)}{$D_p$.2})}~hold.
	\end{proof}

	\begin{theorem}[Error estimate depending on the Sobolev regularity~for~${\Gamma_N=\emptyset}$]\label{error_Sobolev_reg}
		Let $g\hspace*{-0.1em}\in\hspace*{-0.1em} L^\infty(\Omega)$, let $ z\hspace*{-0.1em}\in \hspace*{-0.1em}W^{1,p}(\Omega;\mathbb{R}^d)\cap W^\infty(\textup{div};\Omega)$,  $p\hspace*{-0.1em}\in\hspace*{-0.1em} \left[2,\infty\right)$, with ${\|  z\|_{L^\infty(\Omega;\mathbb{R}^d)}\hspace*{-0.1em}\leq\hspace*{-0.1em} 1}$ be~\mbox{maximal}~for ${D:W^\infty(\textup{div};\Omega)\to\mathbb{R}\cup\{-\infty\}}$, let $u\in \mathcal{A}\hspace*{-0.1em}:=\hspace*{-0.1em}\{{v\hspace*{-0.1em}\in\hspace*{-0.1em} BV(\Omega)\cap L^\infty(\Omega)}\!\mid v\hspace*{-0.1em}=\hspace*{-0.1em}0\text{ in }L^1(\partial\Omega)\}$~be~\mbox{minimal}~for ${I:\mathcal{A}\to \mathbb{R}}$, and let $u_h\in\mathcal{S}^{1,cr}(\mathcal{T}_h)$~be~\mbox{minimal}~for ${I_h:\mathcal{S}^{1,cr}_D(\mathcal{T}_h)\to \mathbb{R}}$. Then, there holds 
		\begin{align*}
			\smash{\|u-\Pi_hu_h\|_{L^2(\Omega)}^2\leq ch^{\frac{p-2}{p-1}}},
		\end{align*}
		where $c>0$ depends only on the quantities $c_{cr}$, $c_d$,~$c_p$,~$\|u\|_{L^\infty(\Omega)}$~and~$\vert \textup{D}u\vert(\Omega)$.
	\end{theorem}
	
	\begin{proof}
		\hspace*{-1mm}Combining the discrete strong coercivity of $I_h\!:\!\mathcal{S}^{1,cr}_D(\mathcal{T}_h)\!\to\! \mathbb{R}$,~i.e.,~\eqref{discrete_strong_coercivity},~and the discrete weak duality principle $I_h(u_h)\hspace*{-0.1em}\ge\hspace*{-0.1em} D_h(\tilde{z}_h)$ for all $\tilde{z}_h\hspace*{-0.1em}\in\hspace*{-0.1em} \mathcal{R}T^0(\mathcal{T}_h)$~(cf.~\eqref{discrete_weak_duality_principle}), we obtain~for~all $\tilde{u}_h\in \mathcal{S}^{1,cr}_D(\mathcal{T}_h)$  and $\tilde{z}_h\in \mathcal{R}T^0(\mathcal{T}_h)$\vspace*{1mm}
		\begin{align}
			\smash{\frac{1}{2}}\|\Pi_h(\tilde{u}_h-u_h)\|_{L^2(\Omega)}^2\leq I_h(\tilde{u}_h)-I_h(u_h)\leq I_h(\tilde{u}_h)-D_h(\tilde{z}_h).\label{eq:errorSobolevreg0}
		\end{align}
		Resorting to Lemma \ref{primal_quasi_interpolant0}, we obtain a function $\tilde{u}_h\hspace*{-0.1em}\in\hspace*{-0.1em}\mathcal{S}^{1,cr}_D(\mathcal{T})$  satisfying \mbox{\textit{(\hyperlink{(P.1)}{P.1})--(\hyperlink{(P.4)}{P.4})}}.
		In addition, Lemma \ref{dual_quasi_interpolant1} yields a vector field $\tilde{z}_h\hspace*{-0.1em}\in\hspace*{-0.1em} \mathcal{R}T^0(\mathcal{T}_h)$ with \textit{(\hyperlink{(Dp.1)}{$D_p$.1})}~and~\textit{(\hyperlink{(Dp.2)}{$D_p$.2})}.
		Combining \textit{(\hyperlink{(P.4)}{P.4})}, \textit{(\hyperlink{(Dp.2)}{$D_p$.2})} and the strong duality principle~${I(u)=D(z)}$~(cf. \eqref{strong_duality_principle}), we deduce from \eqref{eq:errorSobolevreg0} that
		\begin{align}
			\begin{aligned}
			\smash{\frac{1}{2}}\|\Pi_h(\tilde{u}_h-u_h)\|_{L^2(\Omega)}^2&\leq 2c_dc_{cr}\|u\|_{L^\infty(\Omega)}\vert \textup{D}u\vert(\Omega)h+c_ph^{\frac{p-2}{p-1}},
		\end{aligned}\label{eq:errorSobolevreg1}
		\end{align}
		where $c_p\!>\!0$ is as in Lemma \ref{dual_quasi_interpolant1}.
	\!Since ${\tilde{u}_h\!-\!\Pi_h\tilde{u}_h\!=\!\nabla_h \tilde{u}_h\!\cdot\!(\textup{id}_{\mathbb{R}^d}\!-\!\Pi_h\textup{id}_{\mathbb{R}^d})}$~in~$\mathcal{L}^1(\mathcal{T}_h)$, using~\textit{(\hyperlink{(P.1)}{P.1})},~\textit{(\hyperlink{(P.3)}{P.3})} and  $\smash{\|\textup{id}_{\mathbb{R}^d}-\Pi_h\textup{id}_{\mathbb{R}^d}\|_{L^\infty(\Omega;\mathbb{R}^d)}\leq h}$,~we~find~that
		\begin{align}
			\|\tilde{u}_h-\Pi_h\tilde{u}_h\|_{L^2(\Omega)}^2&\leq 2\|\tilde{u}_h\|_{L^\infty(\Omega)}\|\nabla_h \tilde{u}_h\|_{L^1(\Omega;\mathbb{R}^d)}\|\textup{id}_{\mathbb{R}^d}-\Pi_h\textup{id}_{\mathbb{R}^d}\|_{L^\infty(\Omega;\mathbb{R}^d)}\notag\\[-0.25mm]&\leq 2c_d\|u\|_{L^\infty(\Omega)}\vert \textup{D}u\vert(\Omega)h.\label{eq:errorSobolevreg2}
		\end{align}
		Using \textit{(\hyperlink{(P.2)}{P.2})}, \textit{(\hyperlink{(P.3)}{P.3})},  \textit{(\hyperlink{(L0.1)}{L0.1})} and proceeding as for \eqref{eq:errorSobolevreg2}, we further obtain that
		\begin{align}
			\|u-\Pi_h\tilde{u}_h\|_{L^2(\Omega)}^2&\leq \|u-\Pi_h\tilde{u}_h\|_{L^\infty(\Omega)}\big(\|u-\tilde{u}_h\|_{L^1(\Omega)}+\|\tilde{u}_h-\Pi_h\tilde{u}_h\|_{L^1(\Omega)}\big)\notag
			\\[-0.25mm]&\leq (1+c_d)\|u\|_{L^\infty(\Omega)}\big(c_{cr}+1\big)\vert \textup{D}u\vert(\Omega)h.\label{eq:errorSobolevreg3}
		\end{align}
		Finally, combining \eqref{eq:errorSobolevreg1}, \eqref{eq:errorSobolevreg2} and \eqref{eq:errorSobolevreg3}, we conclude~the~claimed~error~bound.
	\end{proof}

	\section{Error estimates for discontinuous dual solutions}\label{sec:non_sobolev}
		
	\qquad In this section, we prove error estimates for the ROF~model~without~\mbox{explicitly} imposing \hspace*{-0.1mm}that \hspace*{-0.1mm}the \hspace*{-0.1mm}dual \hspace*{-0.1mm}solution \hspace*{-0.1mm}possesses \hspace*{-0.1mm}Sobolev \hspace*{-0.1mm}regularity. \hspace*{-0.1mm}Recall~\hspace*{-0.1mm}that,~\hspace*{-0.1mm}in~\hspace*{-0.1mm}\mbox{general}, a solution of the dual ROF model only needs~to~satisfy~${z\!\in\! W^2_N(\textup{div};\Omega)\!\cap\! L^\infty(\Omega;\mathbb{R}^d)}$ with $\|z\|_{L^\infty(\Omega;\mathbb{R}^d)}\!\leq\! 1$. The following lemma gives general assumptions on the dual solution for which it is still possible to construct a~suitable~dual~\mbox{quasi-interpolant}.
	\begin{lemma}[Dual quasi-interpolant for non--Sobolev vector fields]\label{dual_quasi_interpolant}
		Let ${g\in L^2(\Omega)}$ and $z\in W^2_N(\textup{div};\Omega)\cap L^\infty(\Omega;\mathbb{R}^d)$. If $\|\Pi_h I_{\mathcal{R}T}z\|_{L^\infty(\Omega;\mathbb{R}^d)}\leq 1+\kappa(h)$ for $\kappa(h)\ge 0$, then the re-scaled~vector~field $\smash{\tilde{z}_h:=\frac{1}{\gamma_h}I_{\mathcal{R}T}z\in\mathcal{R}T^0_N(\mathcal{T}_h)}$, where $\gamma_h:=1+\kappa(h)>0$, has the~following~properties:
			\begin{description}[noitemsep,topsep=2pt,labelwidth=\widthof{(D.2*)},leftmargin=!,font=\normalfont\itshape]
				\item[(D.1*)] \hypertarget{(D.1*)}{} $\|\Pi_h\tilde{z}_h\|_{L^\infty(\Omega;\mathbb{R}^d)}\leq 1$.
				\item[(D.2*)] \hypertarget{(D.2*)}{} $D_h(\tilde{z}_h)\ge  D( z)-\kappa(h)\|g\|_{L^2(\Omega)}\|\textup{div}( z)\|_{L^2(\Omega)}-\frac{1}{2}\|g-g_h\|_{L^2(\Omega)}^2$. 
			\end{description}
	\end{lemma}

	\begin{proof} Claim \textit{(\hyperlink{(D.1*)}{D.1*})}  is evident. 
		\mbox{Resorting}~to~$\textup{div}(I_{\mathcal{R}T}z)=\Pi_h(\textup{div}(z))$~in~$\mathcal{L}^0(\mathcal{T}_h)$, we deduce that $\smash{\textup{div}(\tilde{z}_h)+g_h=\Pi_h(\frac{1}{\gamma_h}\textup{div}(z)+g)}$ in $\mathcal{L}^0(\mathcal{T}_h)$ and,~hence,~also~using $\|g-g_h\|_{L^2(\Omega)}^2\!=\!\|g\|_{L^2(\Omega)}^2\!-\! \|g_h\|_{L^2(\Omega)}^2$, $I_{K_1(0)}(\Pi_h\tilde{z}_h)\!=\!0$~and~Jensen's~inequality,~that
		\begin{align*}
			D_h(\tilde{z}_h)&= -\tfrac{1}{2}\big\|\Pi_h(\tfrac{1}{\gamma_h}\textup{div}(z)+g)\big\|_{L^2(\Omega)}^2-\tfrac{1}{2}\|g_h\|_{L^2(\Omega)}^2\\[-0.25mm]&
			\ge-\tfrac{1}{2}\big\|\tfrac{1}{\gamma_h}\textup{div}(z)+g\big\|_{L^2(\Omega)}^2+\tfrac{1}{2}\|g\|_{L^2(\Omega)}^2-\tfrac{1}{2}\|g-g_h\|_{L^2(\Omega)}^2\\[-0.25mm]&
			\ge -\tfrac{1}{2}\tfrac{1}{\gamma_h^2}\|\textup{div}(z)\|_{L^2(\Omega)}^2+\tfrac{1}{\gamma_h}(g,\textup{div}(z))_{L^2(\Omega)}-\tfrac{1}{2}\|g-g_h\|_{L^2(\Omega)}^2\\[-0.25mm]&
			\ge D(z)-(1-\tfrac{1}{\gamma_h})(g,\textup{div}(z))_{L^2(\Omega)}-\tfrac{1}{2}\|g-g_h\|_{L^2(\Omega)}^2.
		\end{align*}
	Finally, using that $\smash{\frac{1}{\gamma_h^2}}\leq 1$ and $\smash{1-\frac{1}{\gamma_h}\leq \kappa(h)}$, we conclude that \textit{(\hyperlink{(D.2*)}{D.2*})} holds.\vspace*{-4.5mm}
	\end{proof}\newpage
	
	\begin{theorem}[Error estimate for discontinuous dual solution]\label{non_Sobolev_error_estimate}
		Let $g\in L^\infty(\Omega)$, let $ z\hspace*{-0.15em}\in\hspace*{-0.15em} W^\infty_N(\textup{div};\Omega)$ be maximal for $D\hspace*{-0.15em}:\hspace*{-0.15em}W^\infty_N(\textup{div};\Omega)\hspace*{-0.15em}\to\hspace*{-0.15em} \mathbb{R}\cup\{-\infty\}$ with~the~same~prop-erties as  in Lemma \ref{dual_quasi_interpolant}, let $u\hspace*{-0.15em}\in\hspace*{-0.15em} BV(\Omega)\cap L^\infty(\Omega)$ minimal for ${I\hspace*{-0.15em}:\hspace*{-0.15em}BV(\Omega)\hspace*{-0.15em}\cap \hspace*{-0.15em}L^2(\Omega)\hspace*{-0.15em}\to\hspace*{-0.15em} \mathbb{R}}$, and let $u_h\hspace*{-0.1em}\in\hspace*{-0.1em} \mathcal{S}^{1,cr}(\mathcal{T}_h)$~minimal~for~${I_h\hspace*{-0.1em}:\hspace*{-0.1em}\mathcal{S}^{1,cr}(\mathcal{T}_h)\hspace*{-0.1em}\to\hspace*{-0.1em} \mathbb{R}}$. Then, we have that
		\begin{align*}
			\|u-\Pi_h u_h\|^2_{L^2(\Omega)}\leq c\max\{\kappa(h),h\}.
		\end{align*}
		where $c>0$ depends only on the quantities $c_{cr}$, $c_d$,~$\|u\|_{L^\infty(\Omega)}$,~and~$\vert \textup{D}u\vert(\Omega)$.
	\end{theorem}

	\begin{proof}
		Using the discrete strong coercivity of $I_h:\mathcal{S}^{1,cr}(\mathcal{T}_h)\to \mathbb{R}$,~i.e.,~\eqref{discrete_strong_coercivity}, and the discrete weak duality principle $I_h(u_h)\hspace*{-0.1em}\ge \hspace*{-0.1em}D_h(\tilde{z}_h)$ for all $\tilde{z}_h\hspace*{-0.1em}\in\hspace*{-0.1em} \mathcal{R}T^0_N(\mathcal{T}_h)$~(cf.~\eqref{discrete_weak_duality_principle}), we obtain  for all $\tilde{u}_h\in \mathcal{S}^{1,cr}(\mathcal{T}_h)$  and $\tilde{z}_h\in \mathcal{R}T^0_N(\mathcal{T}_h)$\vspace*{-0.75mm}
		\begin{align}
			\frac{1}{2}\|\Pi_h(\tilde{u}_h-u_h)\|_{L^2(\Omega)}^2\leq I_h(\tilde{u}_h)-I_h(u_h)\leq I_h(\tilde{u}_h)-D_h(\tilde{z}_h).\label{non_Sobolev_error_estimate1}
		\end{align}
		Resorting to Lemma \ref{primal_quasi_interpolant0}, \hspace*{-0.1em}we obtain a function $\tilde{u}_h\hspace*{-0.1em}\in\hspace*{-0.1em}\mathcal{S}^{1,cr}\hspace*{-0.1em}(\mathcal{T})$  satisfying \hspace*{-0.1em}\mbox{\textit{(\hyperlink{(P.1)}{P.1})--(\hyperlink{(P.4)}{P.4})}}.
		In addition, Lemma \ref{dual_quasi_interpolant} yields a vector field $\tilde{z}_h\in \mathcal{R}T^0_N(\mathcal{T}_h)$ with \textit{(\hyperlink{(D.1*)}{D.1*})}~and~\textit{(\hyperlink{(D.2*)}{D.2*}).}
		Then, using \textit{(\hyperlink{(P.4)}{P.4})}, \!\textit{(\hyperlink{(D.2*)}{D.2*})} and the strong duality~principle~${I(u)\!=\!D(z)}$~\!(cf.~\!\eqref{strong_duality_principle}), we deduce from \eqref{non_Sobolev_error_estimate1} that\vspace*{-0.75mm}
		\begin{align*}
			\frac{1}{2}\|\Pi_h(\tilde{u}_h-u_h)\|_{L^2(\Omega)}^2\leq 2c_dc_{cr}\|u\|_{L^\infty(\Omega)}\vert \textup{D}u\vert(\Omega)h+\kappa(h)\|g\|_{L^2(\Omega)}\|\textup{div}( z)\|_{L^2(\Omega)}.
		\end{align*}
		Hence, incorporating  \eqref{eq:errorSobolevreg2}  and \eqref{eq:errorSobolevreg3}, 
		we conclude the claimed error bound.
	\end{proof}

		A sufficient condition for a solution to \eqref{ROF-dual} to guarantee the quasi-optimal~rate $\smash{\mathcal{O}(h^{\frac{1}{2}}\hspace*{-0.1em})}$ is element-wise Lipschitz continuity. In addition, if a solution~to~\eqref{ROF-dual}~is~only element-wise $\alpha$\hspace*{-0.1em}--\hspace*{-0.1em}Hölder continuous, it is, however, possible to derive~the~
		rate~$\smash{\mathcal{O}(h^{\frac{\alpha}{2}})}$.
	
	\begin{lemma}[Dual quasi-interpolant for element-wise $\alpha$--Hölder vector fields]\label{dual_quasi_interpolant_new}
		Let $g\in L^2(\Omega)$ and let $z\in W^2_N(\textup{div};\Omega)\cap L^\infty(\Omega;\mathbb{R}^d)$ be such that $\|z\|_{L^\infty(\Omega;\mathbb{R}^d)}\leq 1$. 
		Furthermore, assume that there exist constants $\alpha\in \left[0,1\right]$  and $c_{\alpha}>0$ such that for all ${T\in \mathcal{T}_h}$, it holds 
		 ${z|_T\in C^{0,\alpha}(T;\mathbb{R}^d)}$ with
		\begin{align}
			\vert z(x)-z(y)\vert \leq c_{\alpha}\vert x-y\vert^{\alpha}\label{hoelder}
		\end{align}
		for \hspace*{-0.1mm}all \hspace*{-0.1mm}$x,y\!\in\! T$. \hspace*{-0.3em}Then,  \hspace*{-0.1mm}the \hspace*{-0.1mm}assumptions \hspace*{-0.1mm}in  \hspace*{-0.1mm}Lemma \hspace*{-0.1mm}\ref{dual_quasi_interpolant}~\hspace*{-0.1mm}are~\hspace*{-0.1mm}satisfied~\hspace*{-0.1mm}with~\hspace*{-0.1mm}$\smash{\kappa(h)\!=\!\mathcal{O}(h^{\alpha})}$.
	\end{lemma}

	\begin{remark}\label{resolved_discontinuity}
			Lemma \ref{dual_quasi_interpolant_new} is of particular interest if the discontinuity set $J_z$ of a piece-wise regular (piece-wise Lipschitz or piece-wise $\alpha$--Hölder~continuous) vector field $z\!\in\! W^2_N(\textup{div};\Omega)\cap L^\infty(\Omega;\mathbb{R}^d)$  is resolved~by~the~\mbox{triangulation},~i.e.,~${J_z\!\subseteq\! \bigcup_{S\in \mathcal{S}_h}{\!S}}$.
	\end{remark}

	\begin{proof} (of Lemma \ref{dual_quasi_interpolant_new})
		We  need to check that $\|\Pi_h I_{\mathcal{R}T}z\|_{L^\infty(\Omega;\mathbb{R}^d)}\!\leq\! 1+\kappa(h)$~for some \hspace*{-0.1mm}${\kappa(h)\hspace*{-0.2em}>\hspace*{-0.2em}-1}$~\hspace*{-0.1mm}with~\hspace*{-0.1mm}${\kappa(h)\hspace*{-0.2em}=\hspace*{-0.2em}\mathcal{O}(h^{\alpha})}$.
		\hspace*{-0.5mm}Note \hspace*{-0.1mm}that \hspace*{-0.1mm}$I_{\mathcal{R}T}(z(x_T))\hspace*{-0.2em}=\hspace*{-0.2em}z(x_T)$~\hspace*{-0.1mm}in~\hspace*{-0.1mm}$T$~\hspace*{-0.1mm}for~\hspace*{-0.1mm}all~\hspace*{-0.1mm}${T\hspace*{-0.2em}\in \hspace*{-0.2em}\mathcal{T}_h}$, which  results~from 
		$\textup{div}(I_{\mathcal{R}T}(z(x_T)))=\Pi_h(\textup{div}(z(x_T)))=0$ in $T$ for all $T\in \mathcal{T}_h$.
Using~this, \textit{(\hyperlink{RT.2}{RT.2})}, \eqref{hoelder} and that $\|z\|_{L^\infty(\Omega;\mathbb{R}^d)}\leq 1$, we deduce that for all $T\in \mathcal{T}_h$
		\begin{align}
			\begin{aligned}
			\vert  (I_{\mathcal{R}T}z)(x_T)\vert&\leq \vert  I_{\mathcal{R}T}(z-z(x_T))(x_T)\vert +\vert z(x_T)\vert\\&\leq
			 \| I_{\mathcal{R}T}(z-z(x_T))\|_{L^\infty(T;\mathbb{R}^d)}+1\\&\leq
			 c_{\mathcal{R}T}\|z-z(x_T)\|_{L^\infty(T;\mathbb{R}^d)}+1
			 \\&\leq
			 c_{\mathcal{R}T}{\sup}_{x\in T}{\vert x-x_T\vert^\alpha} +1
			 \\&\leq
			 c_{\mathcal{R}T}c_\alpha h_T^{\alpha} +1,
\end{aligned}\label{dual_quasi_interpolant_new.2}
		\end{align}
		i.e., setting $\kappa(h)\hspace*{-0.1em}:=\hspace*{-0.1em}c_{\mathcal{R}T}c_\alpha h^\alpha$,  we conclude that $\|\Pi_h I_{\mathcal{R}T}z\|_{L^\infty(\Omega;\mathbb{R}^d)}\hspace*{-0.1em}\leq\hspace*{-0.1em} 1+\kappa(h)$.
\end{proof}\newpage

	\begin{theorem}[Error estimate for element-wise $\alpha$--Hölder  dual solution]\label{error_elementwise_hoelder_reg}
		Let $z\in W^2_N(\textup{div};\Omega)\cap L^\infty(\Omega;\mathbb{R}^d)$  be maximal for $D:W^2_N(\textup{div};\Omega)\cap L^\infty(\Omega;\mathbb{R}^d)\to \mathbb{R}\cup\{-\infty\}$ with the same \mbox{properties}  as  in Lemma \ref{dual_quasi_interpolant_new}, let $u\in BV(\Omega)\cap L^\infty(\Omega)$ minimal for ${I\hspace*{-0.1em}:\hspace*{-0.1em}BV(\Omega)\cap L^2(\Omega)\hspace*{-0.1em}\to\hspace*{-0.1em} \mathbb{R}}$ and let $u_h\hspace*{-0.15em}\in\hspace*{-0.15em} \mathcal{S}^{1,cr}(\mathcal{T}_h)$~minimal~for~${I_h\hspace*{-0.15em}:\hspace*{-0.15em}\mathcal{S}^{1,cr}(\mathcal{T}_h)\hspace*{-0.15em}\to\hspace*{-0.15em} \mathbb{R}}$. Then, we have that\vspace*{-0.5mm}
		\begin{align*}
			\|u-\Pi_h u_h\|^2_{L^2(\Omega)}\leq ch^{\alpha}.
		\end{align*}
		where $c\!>\!0$ depends only on the quantities $c_{cr}$, $c_d$, $c_\alpha$,~$\|u\|_{L^\infty(\Omega)}$,~and~$\vert \textup{D}u\vert(\Omega)$.
	\end{theorem}

	\begin{proof}
		Follows from Theorem \ref{non_Sobolev_error_estimate} by resorting to Lemma \ref{dual_quasi_interpolant_new}.
	\end{proof}

	\begin{remark}[Comparison of Theorem \ref{error_Sobolev_reg} and Theorem \ref{error_elementwise_hoelder_reg}]
		\begin{description}[noitemsep,topsep=2pt,labelwidth=\widthof{(ii)},leftmargin=!,font=\normalfont\itshape]
			\item[(i)] \hspace*{-2mm}If $\alpha\hspace*{-0.1em}=\hspace*{-0.1em}1$, then Theorem \ref{error_elementwise_hoelder_reg} extends the results \cite[Proposition 4.2]{Bar21} and \cite[Sectioin 5.1.1]{CP20} 
			to the case~of~an~existing element-wise Lipschitz continuous solution to \eqref{ROF-dual}.
			
			\item[(ii)] If $p>d$ in Theorem \ref{error_Sobolev_reg}, then $z\in W^{1,p}(\Omega;\mathbb{R}^d)$ satisfies $z\in C^{0,\alpha}(\overline{\Omega};\mathbb{R}^d)$ for $\alpha=1-\frac{d}{p}$ by Sobolev's embedding theorem \cite[Corollary~9.14]{Bre10}.~As~a~result, 
			Theorem \hspace*{-0.1mm}\ref{error_elementwise_hoelder_reg} \hspace*{-0.1mm}in \hspace*{-0.1mm}the \hspace*{-0.1mm}particular \hspace*{-0.1mm}case \hspace*{-0.1mm}$\Gamma_{\!N}\!=\!\emptyset$ \hspace*{-0.1mm}is \hspace*{-0.1mm}applicable~\hspace*{-0.1mm}and~\hspace*{-0.1mm}yields~\hspace*{-0.1mm}the~\hspace*{-0.1mm}rate~\hspace*{-0.1mm}$\smash{\mathcal{O}(\hspace*{-0.1em}h^{\hspace*{-0.1em}\frac{\alpha}{2}}\hspace*{-0.1em})}$. On the other hand, Theorem \ref{error_Sobolev_reg} yields the slightly improved  rate  $\smash{\mathcal{O}(h^{\smash{\frac{p-2}{2(p-1)}}})}$, which gives the impression that Theorem \ref{error_Sobolev_reg} is utterly superior~to~\mbox{Theorem~\!\ref{error_elementwise_hoelder_reg}}.  Nevertheless, the major strength of Theorem~\ref{error_elementwise_hoelder_reg} -- and equally of  Lemma~\ref{dual_quasi_interpolant_new}~-- is that it is also applicable when it is unclear whether a solution 
			to \eqref{ROF-dual}~with Sobolev regularity is available. This allows us to justify analytically the quasi-optimal rate $\smash{\mathcal{O}(h^{\frac{1}{2}})}$ for the setting in Section \ref{subsec:irregular_example} at least for the particular case that the discontinuity set $J_z$ is resolved~by~the~\mbox{triangulation},~i.e.,~${J_z\!\subseteq\! \bigcup_{S\in \mathcal{S}_h}{\!S}}$, cf. Example \ref{two_disk_problem} and Example \ref{four_disk_problem}.
		\end{description}
	\end{remark}

	If the discontinuity set of a solution to \eqref{ROF-dual} is not resolved by the triangulation, then, apparently, Theorem \ref{error_elementwise_hoelder_reg} does not apply. In this case, however, the following 
	argument applies, which exploits that~for~${g\hspace*{-0.1em}\in L^\infty(\Omega)}$,~we~have~that~${\textup{div}(z)\hspace*{-0.1em}\in\hspace*{-0.1em} L^\infty(\Omega)}$, which to some extent can serve  as a substitute~for~${\nabla z\in L^\infty(\Omega;\mathbb{R}^{d\times d})}$.
	
	\begin{remark}[Optimal dual quasi-interpolant for non--Lipschitz vector fields]\label{soeren_trick}
		Let $z\in W^\infty_N(\textup{div};\Omega)$ be such that $\|z\|_{L^\infty(\Omega;\mathbb{R}^d)}\leq 1$. Furthermore, assume that there exists a constant $\tilde{c}_z>0$ such that
		for  all $T\in \mathcal{T}_h$,  
		there exists some $\widetilde x_T\in T$ such that
		\begin{align}
			\vert (I_{\mathcal{R}T}z)(\widetilde x_T)\vert \leq 1+\tilde{c}_zh,\label{eq:soeren_trick}
		\end{align}
		For each $T\in \mathcal{T}_h$, since $I_{\mathcal{R}T}z\in \mathcal{R}T^0_N(\mathcal{T}_h)\subseteq \mathcal{L}^1(\mathcal{T}_h)^d$, we have that
		\begin{align}
			(I_{\mathcal{R}T}z)(x)=(I_{\mathcal{R}T}z)(x_T)+d^{-1}\textup{div}(I_{\mathcal{R}T}z)(x-x_T)\label{eq:soeren_trick1}
		\end{align}
		for \hspace*{-0.1mm}all \hspace*{-0.1mm}$x\!\in\! T$. \hspace*{-0.1mm}Thus, \hspace*{-0.1mm}resorting \hspace*{-0.1mm}to  \hspace*{-0.1mm}${\textup{div}(I_{\mathcal{R}T}z)\!=\!\Pi_h(\textup{div}(z))}$~\hspace*{-0.1mm}in~\hspace*{-0.1mm}$\mathcal{L}^0(\mathcal{T}_h)$,~\hspace*{-0.1mm}also~\hspace*{-0.1mm}\mbox{using}~\hspace*{-0.1mm}\eqref{eq:soeren_trick} and \textit{(\hyperlink{(L0.1)}{L0.1})} in \eqref{eq:soeren_trick1} at $x=\widetilde x_T\in T$,~we~conclude~that
		\begin{align*}
			\|\Pi_hI_{\mathcal{R}T}z\|_{L^\infty(T;\mathbb{R}^d)}&\leq \vert (I_{\mathcal{R}T}z)(\tilde x_T)\vert +d^{-1}\|\Pi_h(\textup{div}(z))\|_{L^\infty(T)}\vert \tilde x_T-x_T\vert \\&\leq
			1+\tilde{c}_zh+d^{-1}\|\textup{div}(z)\|_{L^\infty(T)} h_T,
		\end{align*}
		i.e., we have that $\|\Pi_hI_{\mathcal{R}T}z\|_{L^\infty(T;\mathbb{R}^d)}\leq 1+\big(\tilde{c}_z+d^{-1}\|\textup{div}(z)\|_{L^\infty(\Omega)}\big)h$.
\end{remark}
	
	The following remark discusses particular sufficient conditions  for \eqref{eq:soeren_trick} on a vector field $z\in W^\infty_N(\textup{div};\Omega)$ that is piece-wise Lipschitz continuous, such as, e.g., that
	its discontinuity set $J_z$ is approximated by $\mathcal{T}_h$, $h>0$, with rate $\mathcal{O}(h)$ or that
	$\vert z\vert <1$ along $J_z$ while, simultaneously, its jump $[\![z]\!]$ over $J_z$ remains small.~~~~~~~~~~~\newpage
	On the other hand, this remark finds that \eqref{eq:soeren_trick}   cannot be expected, in general, for piece-wise Lipschitz continuous vector fields, even in generic situations.~~~~~~~~~~~~
	 
	\begin{remark}[Sufficient conditions for \eqref{eq:soeren_trick}]\label{rem:soeren_trick_special_cases}
		Let $d=2$ and $z\in W^\infty_N(\textup{div};\Omega)$~with $\|z\|_{L^\infty(\Omega;\mathbb{R}^2)}\leq 1$ be piece-wise Lipschitz continuous, i.e., there~exist~open~$\Omega_i\subseteq \Omega$, $i\!=\!1,\dots,m$, $m\!\in\! \mathbb{N}$, with
		$z|_{\Omega_i}\!\in\! W^{1,\infty}(\Omega_i;\mathbb{R}^2)$~for~all~$i\!=\!1,\dots,m$~and~${\overline{\Omega}\!=\!\bigcup_{i=1}^m{\overline{\Omega_i}}}$. 
		Next, we fix an arbitrary $T\in \mathcal{T}_h$.
		Then, we need to distinguish two cases:
		\begin{description}[noitemsep,topsep=2pt,labelwidth=\widthof{(ii)},leftmargin=!,font=\normalfont\itshape]
				\item[(i)] Assume that $T\subseteq \overline{\Omega_i}$ for some $i=1,\dots,m$.  Then, we~deduce~along~the~lines~of the proof of \eqref{dual_quasi_interpolant_new.2} that
				\begin{align*}
					\vert (I_{\mathcal{R}T}z)(x_T)\vert \leq 1+ c_{\mathcal{R}T} \|\nabla z\|_{L^\infty(\Omega_i;\mathbb{R}^{2\times 2})} h_T,
				\end{align*}
				i.e., \hspace*{-0.1em}$\vert (I_{\mathcal{R}T}z)(x_T)\vert \!\leq\!  1+c_zc_{\mathcal{R}T}h_T$, where $c_z\!:=\!\max_{i=1,\dots,m}{\!\|\nabla (\hspace*{-0.05em}z|_{\Omega_i}\hspace*{-0.05em})\|_{L^\infty(\Omega;\mathbb{R}^{2\times 2})}}$.
				\item[(ii)] Assume that there exists an interface $\gamma=\partial\Omega_a\cap \partial\Omega_b$~for~some~${a,b=1,\dots,m}$ such that $\textup{int}(T)\cap \gamma\neq \emptyset$\footnote{Apparently, we should also take into account the case in which $T\in \mathcal{T}_h$ is intersected by two or more interfaces. However, 
					for the benefit of readability, we limit ourselves to this simplified case.\vspace*{-7mm}}. As $z|_{\smash{\Omega_a}}\in W^{1,\infty}(\Omega_a;\mathbb{R}^2)$ and $z|_{\smash{\Omega_b}}\in W^{1,\infty}(\Omega_b;\mathbb{R}^2)$, without loss of generality, we may assume that~${\gamma \subseteq b_\gamma+\mathbb{R}t_\gamma}$~for~some  $b_\gamma\in \mathbb{R}^2$ and ${t_\gamma\in \mathbb{S}^1}$. Next, fix  $x_\gamma\in \gamma $ and set ${z_a:=(z|_{\smash{\overline{\Omega_a}}})(x_\gamma), z_b:=(z|_{\smash{\overline{\Omega_b}}})(x_\gamma)\in \mathbb{R}^2}$. Then, for $i\in \{a,b\}$, it holds
				\begin{align*}
					\vert z_i-(z|_{\smash{\overline{\Omega_i}}})(x)\vert& \leq \|\nabla (z|_{\smash{\Omega_i}})\|_{L^\infty(\Omega_i;\mathbb{R}^{2\times 2})}\vert x_{\gamma}-x\vert\leq c_zh_T\;\;\text{ for all }x\in T\cap \overline{\Omega_i}.
				\end{align*}
				Furthermore, if $n_\gamma\in \mathbb{S}^1$ denotes a unit normal to $t_\gamma\in \mathbb{S}^1$, i.e., $n_\gamma \cdot t_\gamma=0$, then, taking into account that $z\in W^\infty_N(\textup{div};\Omega)$,~we find~that~${z_a\cdot n_\gamma=z_b\cdot n_\gamma}$. 
				Thus, if we define $z_T(x):=z_a$ for $x\in T\cap \overline{\Omega_a}$ and $z_T(x):=z_b$ for $x\in T\cap \overline{\Omega_b}$, cf. Figure \ref{fig:triangles}, \textit{($\alpha$)}, then 
				$z_T\in W^{\infty}(\textup{div};T)$ and, owing to \textit{(\hyperlink{(RT.2)}{RT.2})}, 
				\begin{align*}
					\begin{aligned}
						\|I_{\mathcal{R}T}z_T-I_{\mathcal{R}T}z\|_{L^\infty(T;\mathbb{R}^2)}
						\leq c_{\mathcal{R}T}\|z_T-z\|_{L^\infty(T;\mathbb{R}^2)}\leq c_{\mathcal{R}T}c_zh_T,
					\end{aligned}
				\end{align*}
				i.e., we have that
				\begin{align}
					\|I_{\mathcal{R}T}z\|_{L^\infty(T;\mathbb{R}^2)}\leq \|I_{\mathcal{R}T}z_T\|_{L^\infty(T;\mathbb{R}^2)}+c_{\mathcal{R}T}c_zh_T.\label{eq:soeren_trick_special_cases.1}
				\end{align}
				As a result of \eqref{eq:soeren_trick_special_cases.1}, it is sufficient to prove ${\|I_{\mathcal{R}T}z_T\|_{L^\infty(T;\mathbb{R}^2)}\leq 1+\mathcal{O}(h)}$ to conclude that $\|I_{\mathcal{R}T}z\|_{L^\infty(T;\mathbb{R}^2)}\leq 1+\mathcal{O}(h)$.
				Because, owing to Lemma~\ref{dual_quasi_interpolant}~(ii), it holds
				\begin{align*}
					\textup{div}(I_{\mathcal{R}T}z_T)=\Pi_h(\textup{div}(z_T))=0\quad\text{ in }T,
				\end{align*}
				where $I_{\mathcal{R}T}z_T:=\sum_{S\in \mathcal{S}_h;S\subseteq \partial T}{z_T\cdot n_S\psi_S}$, it even holds ${I_{\mathcal{R}T}z_T\equiv\textup{const}}$~in~$T$.
				Next, we denote by $S_1\in \mathcal{S}_h$ a side of $T\in \mathcal{T}_h$ such that $S_1\cap \gamma \neq \emptyset$ and by $S_2\in \mathcal{S}_h$ the side~of~${T\in \mathcal{T}_h}$~such~that~${S_2\cap \gamma =  \emptyset}$. Let $n_1,n_2\in \mathbb{S}^1$ denote the corresponding unit normal~vectors~to~${S_1,S_2\in \mathcal{T}_h}$, resp., cf. Figure \ref{fig:triangles}, \textit{($\alpha$)}.  Then, it holds
				\begin{align}
					\left.\begin{aligned}
						I_{\mathcal{R}T}z_T\cdot n_1& = \int_{S_1}{z_T\cdot n_1\,\textup{d}s}=\frac{\vert S_1\cap \Omega_b\vert }{\vert S_1\vert}z_b\cdot n_1+\frac{\vert S_1\cap \Omega_a\vert }{\vert S_1\vert}z_a\cdot n_1,\\
					I_{\mathcal{R}T}z_T\cdot n_2& = \int_{S_2}{z_T\cdot n_2\,\textup{d}s}=z_b\cdot n_2.
				\end{aligned}\;\right\} \label{eq:soeren_trick_special_cases.2}
				\end{align}
				Introducing $\rho:=\vert S_1\cap \Omega_b\vert /\vert S_1\vert\in \left[0,1\right]$ as well as ${M_T:=(n_1,n_2)\in \mathbb{R}^{2\times 2}}$, also exploiting that $z_a=z_b+((z_a-z_b)\cdot t_\gamma)t_\gamma$, 
				where we used that ${z_a\cdot n_\gamma=z_b\cdot n_\gamma}$, the system \eqref{eq:soeren_trick_special_cases.2} can be rewritten  as 
				\begin{align*}
					M^\top_T I_{\mathcal{R}T}z_T= M^\top_T z_b+(1-\rho)((z_a-z_b)\cdot t_\gamma)(t_\gamma \cdot n_1)e_1,
				\end{align*}
				i.e., since $M_T^\top\in \mathbb{R}^{2\times 2}$ is a regular matrix, we find that 
				\begin{align}
					I_{\mathcal{R}T}z_T= z_b+(1-\rho)((z_a-z_b)\cdot t_\gamma)(t_\gamma \cdot n_1)M_T^{-\top}e_1.\label{eq:soeren_trick_special_cases.3}
				\end{align}
				Resorting to the formula \eqref{eq:soeren_trick_special_cases.3}, we can derive special~cases~that~imply~\eqref{eq:soeren_trick}:
				\begin{description}[noitemsep,topsep=2pt,labelwidth=\widthof{(ii.a)},leftmargin=!,font=\normalfont\itshape]
					\item[(ii.a)]  If $t_\gamma \cdot n_1=\mathcal{O}(h)$, i.e., $(b_\gamma+\mathbb{R}t_\gamma)\cap T$ approximates $S_1$ with rate~$\mathcal{O}(h)$, cf.~Figure~\ref{fig:triangles}, \textit{($\beta$)}, then
					$	\vert I_{\mathcal{R}T}z_T\vert \leq  \vert z_b\vert +\mathcal{O}(h)\leq 1+\mathcal{O}(h)$.
					
					\item[(ii.b)] If $1-\rho=\mathcal{O}(h)$, i.e.,  $S_1\cap \Omega_b$ approximates $S_1$ with rate~$\mathcal{O}(h)$, cf. Figure \ref{fig:triangles}, \textit{($\gamma$)}, then
					$	\vert I_{\mathcal{R}T}z_T\vert \leq  \vert z_b\vert +\mathcal{O}(h)\leq 1+\mathcal{O}(h)$.
				\end{description}
				Apparently,  \textit{(ii.a)} and \textit{(ii.b)} describe the particular case in which the disconti-nuity \hspace*{-0.1mm}set \hspace*{-0.1mm}$J_z$ \hspace*{-0.1mm}is \hspace*{-0.1mm}not \hspace*{-0.1mm}resolved \hspace*{-0.1mm}by \hspace*{-0.1mm}the \hspace*{-0.1mm}triangulation \hspace*{-0.1mm}but~\hspace*{-0.1mm}approximated~\hspace*{-0.1mm}with~\hspace*{-0.1mm}rate~\hspace*{-0.1mm}$\mathcal{O}(h)$.
				\begin{description}[noitemsep,topsep=2pt,labelwidth=\widthof{(ii.a)},leftmargin=!,font=\normalfont\itshape]
					\item[(ii.c)] If we have that both $\vert z_b\vert<1 $ and  $(z_a-z_b)\cdot t_\gamma$ is~sufficiently~small, i.e., such that $	\vert (1-\rho)((z_a-z_b)\cdot t_\gamma)M^{-\top}_T(t_\gamma \cdot n_1)e_1\vert \leq  1-\vert z_b\vert+\mathcal{O}(h)$,
					then  $	\vert I_{\mathcal{R}T}z_T\vert \leq  \vert z_b\vert +1-\vert z_b\vert +\mathcal{O}(h)= 1+\mathcal{O}(h)$.
					\item[(ii.d)] If $T$ is nearly right-angled, so that $M_T$ is approximately an orthogonal matrix, i.e., $ M_T^{-\top} = M_T+\mathcal{O}(h)$, 
					and ${t_\gamma =\pm n_1+\mathcal{O}(h)}$,  then, using that $z_a\cdot n_2=z_b\cdot n_2+\mathcal{O}(h)$ because $n_\gamma = \pm n_2+\mathcal{O}(h)$,~we~deduce~that $z_b=(z_b\cdot n_1)n_1+(1-\rho)(z_a\cdot n_2)n_2+\rho(z_b\cdot n_2)n_2+\mathcal{O}(h)$~and,~thus,~that
					\begin{align*}
					I_{\mathcal{R}T}z_T&=z_b+(1-\rho)((z_a-z_b)\cdot n_1)n_1+\mathcal{O}(h)\\&
					=\rho((z_b\cdot n_1)n_1+\rho(z_b\cdot n_2)n_2)\\&\quad+(1-\rho)((z_a\cdot n_1)n_1+(z_a\cdot n_2)n_2)+\mathcal{O}(h)\\&
					=(1-\rho)z_a+\rho z_b+\mathcal{O}(h),
					\end{align*}
					which implies that $\vert I_{\mathcal{R}T}z_T\vert \leq (1-\rho)\vert z_a\vert +\rho \vert z_b\vert+\mathcal{O}(h)\leq 1+\mathcal{O}(h)$.
				\end{description}
				More generally, the sub-cases \textit{(ii.a)--(ii.d)} can occur in combination so that the conclusion \hspace*{-0.1mm}holds \hspace*{-0.1mm}under \hspace*{-0.1mm}significantly \hspace*{-0.1mm}weaker \hspace*{-0.1mm}conditions \hspace*{-0.1mm}on \hspace*{-0.1mm}the~\hspace*{-0.1mm}\mbox{individual}~\hspace*{-0.1mm}\mbox{factors}.
				On the other hand,  the formula \eqref{eq:soeren_trick_special_cases.3}, simultaneously, demonstrates that $ \|I_{\mathcal{R}T}z_T\|_{L^\infty(T;\mathbb{R}^2)}\leq 1+\mathcal{O}(h) $ and, therefore, also ${\|I_{\mathcal{R}T}z\|_{L^\infty(T;\mathbb{R}^2)}\leq 1+\mathcal{O}(h)}$ cannot be expected in general, even in generic situations.
		\end{description}
	\begin{figure}[h]
	\centering

\tikzset{every picture/.style={line width=0.75pt}} 
\vspace*{-0.35cm}
\begin{tikzpicture}[x=0.8pt,y=0.8pt,yscale=-1,xscale=1]
	
	\draw [color={rgb, 255:red, 155; green, 155; blue, 155 }  ,draw opacity=1 ]   (127.41,87) -- (95.41,139) ;
	\draw [color={rgb, 255:red, 208; green, 2; blue, 27 }  ,draw opacity=1 ][line width=0.75]    (181.41,1) -- (127.41,87) -- (51.41,211) ;
	\draw  [color={rgb, 255:red, 155; green, 155; blue, 155 }  ,draw opacity=1 ] (95.74,138.33) .. controls (95.71,143) and (98.03,145.35) .. (102.69,145.38) -- (143.63,145.65) .. controls (150.29,145.7) and (153.61,148.05) .. (153.58,152.72) .. controls (153.61,148.05) and (156.95,145.74) .. (163.62,145.78)(160.62,145.76) -- (205.18,146.06) .. controls (209.85,146.09) and (212.2,143.78) .. (212.23,139.11) ;
	\draw    (35.41,87) -- (11.83,96.27) ;
	\draw [shift={(9.96,97)}, rotate = 338.54] [fill={rgb, 255:red, 0; green, 0; blue, 0 }  ][line width=0.08]  [draw opacity=0] (9.6,-2.4) -- (0,0) -- (9.6,2.4) -- cycle    ;
	\draw    (106.23,139) -- (106.23,165) ;
	\draw [shift={(106.23,167)}, rotate = 270] [fill={rgb, 255:red, 0; green, 0; blue, 0 }  ][line width=0.08]  [draw opacity=0] (9.6,-2.4) -- (0,0) -- (9.6,2.4) -- cycle    ;
	\draw   (215.96,14.66) -- (375.18,92.37) -- (251.23,92.37) -- cycle ;
	\draw [color={rgb, 255:red, 208; green, 2; blue, 27 }  ,draw opacity=1 ]   (225.27,111) -- (369.27,201) ;
	\draw   (217.62,115.82) -- (376.83,193.53) -- (252.89,193.53) -- cycle ;
	\draw  [color={rgb, 255:red, 155; green, 155; blue, 155 }  ,draw opacity=1 ] (359.27,205.67) .. controls (359.27,208.14) and (360.51,209.37) .. (362.98,209.37) -- (362.98,209.37) .. controls (366.51,209.37) and (368.27,210.61) .. (368.27,213.08) .. controls (368.27,210.61) and (370.04,209.37) .. (373.57,209.37)(371.98,209.37) -- (373.57,209.37) .. controls (376.04,209.37) and (377.27,208.14) .. (377.27,205.67) ;
	\draw [color={rgb, 255:red, 208; green, 2; blue, 27 }  ,draw opacity=1 ]   (393.62,69) -- (233.62,103) ;
	\draw    (111.41,113) -- (89.72,98.64) ;
	\draw [shift={(88.05,97.54)}, rotate = 33.51] [fill={rgb, 255:red, 0; green, 0; blue, 0 }  ][line width=0.08]  [draw opacity=0] (9.6,-2.4) -- (0,0) -- (9.6,2.4) -- cycle    ;
	\draw    (111.41,113) -- (126.36,88.7) ;
	\draw [shift={(127.41,87)}, rotate = 121.61] [fill={rgb, 255:red, 0; green, 0; blue, 0 }  ][line width=0.08]  [draw opacity=0] (9.6,-2.4) -- (0,0) -- (9.6,2.4) -- cycle    ;
	\draw [color={rgb, 255:red, 155; green, 155; blue, 155 }  ,draw opacity=1 ] [dash pattern={on 0.84pt off 2.51pt}]  (301.67,104.67) -- (286.01,96.72) ;
	\draw  [draw opacity=0] (291.59,90.35) .. controls (291.91,91.57) and (292.01,92.87) .. (291.84,94.19) .. controls (291.18,99.29) and (286.76,103) .. (281.66,103) -- (281.62,92.87) -- cycle ; \draw   (291.59,90.35) .. controls (291.91,91.57) and (292.01,92.87) .. (291.84,94.19) .. controls (291.18,99.29) and (286.76,103) .. (281.66,103) ;
	\draw    (281.62,92.87) -- (281.69,112.99) ;
	\draw [shift={(281.7,114.99)}, rotate = 269.8] [fill={rgb, 255:red, 0; green, 0; blue, 0 }  ][line width=0.08]  [draw opacity=0] (9.6,-2.4) -- (0,0) -- (9.6,2.4) -- cycle    ;
	\draw    (281.62,92.87) -- (305.23,87.4) ;
	\draw [shift={(307.18,86.95)}, rotate = 166.97] [fill={rgb, 255:red, 0; green, 0; blue, 0 }  ][line width=0.08]  [draw opacity=0] (9.6,-2.4) -- (0,0) -- (9.6,2.4) -- cycle    ;
	\draw [color={rgb, 255:red, 155; green, 155; blue, 155 }  ,draw opacity=1 ] [dash pattern={on 0.84pt off 2.51pt}]  (377.27,205) -- (376.83,193.53) ;
	\draw [color={rgb, 255:red, 155; green, 155; blue, 155 }  ,draw opacity=1 ] [dash pattern={on 0.84pt off 2.51pt}]  (359.27,205) -- (359.27,194.87) ;
	\draw  [color={rgb, 255:red, 155; green, 155; blue, 155 }  ,draw opacity=1 ] (54.74,138.33) .. controls (54.67,143) and (56.96,145.37) .. (61.63,145.45) -- (65.01,145.5) .. controls (71.68,145.61) and (74.97,148) .. (74.9,152.67) .. controls (74.97,148) and (78.34,145.72) .. (85.01,145.83)(82.01,145.78) -- (88.45,145.89) .. controls (93.12,145.96) and (95.49,143.67) .. (95.56,139) ;
	\draw   (9.96,19) -- (210.97,139) -- (54.5,139) -- cycle ;
	\draw [color={rgb, 255:red, 74; green, 74; blue, 74 }  ,draw opacity=1 ]   (139.71,19.07) -- (139.11,49.91) ;
	\draw [shift={(139.07,51.91)}, rotate = 271.12] [fill={rgb, 255:red, 74; green, 74; blue, 74 }  ,fill opacity=1 ][line width=0.08]  [draw opacity=0] (9.6,-2.4) -- (0,0) -- (9.6,2.4) -- cycle    ;
	\draw [color={rgb, 255:red, 74; green, 74; blue, 74 }  ,draw opacity=1 ]   (152.71,61.07) -- (180.71,60.14) ;
	\draw [shift={(182.71,60.07)}, rotate = 178.09] [fill={rgb, 255:red, 74; green, 74; blue, 74 }  ,fill opacity=1 ][line width=0.08]  [draw opacity=0] (9.6,-2.4) -- (0,0) -- (9.6,2.4) -- cycle    ;
	
	\draw (67.94,77.27) node [anchor=north west][inner sep=0.75pt]    {$T$};
	\draw (120,30) node [anchor=north west][inner sep=0.75pt]  [font=\normalsize]  {$z_{b}$};
	\draw (160,68) node [anchor=north west][inner sep=0.75pt]  [font=\normalsize]  {$z_{a}$};
	\draw (35,175) node [anchor=north west][inner sep=0.75pt]  [font=\normalsize]  {$\Omega _{b}$};
	\draw (85,175) node [anchor=north west][inner sep=0.75pt]  [font=\normalsize]  {$\Omega _{a}$};
	\draw (112,110) node [anchor=north west][inner sep=0.75pt]  [font=\normalsize]  {$x_{\gamma }$};
	\draw (69.23,153) node [anchor=north west][inner sep=0.75pt]  [font=\normalsize]  {$\rho $};
	\draw (138.5,151.4) node [anchor=north west][inner sep=0.75pt]  [font=\normalsize]  {$1-\rho $};
	\draw (115.23,123) node [anchor=north west][inner sep=0.75pt]  [font=\normalsize]  {$S_{1}$};
	\draw (37.23,76) node [anchor=north west][inner sep=0.75pt]  [font=\normalsize]  {$S_{2}$};
	\draw (13.23,100) node [anchor=north west][inner sep=0.75pt]  [font=\normalsize]  {$n_{2}$};
	\draw (109,160) node [anchor=north west][inner sep=0.75pt]  [font=\normalsize]  {$n_{1}$};
	\draw (22.5,210) node [anchor=north west][inner sep=0.75pt]  [font=\normalsize]  {$J_{z_T}=\mathbb{R} t_\gamma $};
	\draw (110.23,203.4) node [anchor=north west][inner sep=0.75pt]    {$(\alpha)$};
	\draw (260.62,78) node [anchor=north west][inner sep=0.75pt]  [font=\normalsize]  {$S_{1}$};
	\draw (240,48.4) node [anchor=north west][inner sep=0.75pt]  [font=\normalsize]  {$S_{2}$};
	\draw (331.5,215.4) node [anchor=north west][inner sep=0.75pt]  [font=\normalsize]  {$1-\rho =\mathcal{O}(h)$};
	\draw (299.41,178) node [anchor=north west][inner sep=0.75pt]  [font=\normalsize]  {$S_{1}$};
	\draw (240,150.02) node [anchor=north west][inner sep=0.75pt]  [font=\normalsize]  {$S_{2}$};
	\draw (210,95) node [anchor=north west][inner sep=0.75pt]  [font=\normalsize]  {$\mathbb{R} t_\gamma $};
	\draw (83.23,104) node [anchor=north west][inner sep=0.75pt]  [font=\normalsize]  {$n_{\gamma }$};
	\draw (126.23,96.4) node [anchor=north west][inner sep=0.75pt]  [font=\normalsize]  {$t_{\gamma }$};
	\draw (340,30) node [anchor=north west][inner sep=0.75pt]    {$( \beta)$};
	\draw (340,130) node [anchor=north west][inner sep=0.75pt]    {$( \gamma)$};
	\draw (265,108.4) node [anchor=north west][inner sep=0.75pt]  [font=\normalsize]  {$n_{1}$};
	\draw (297.62,72) node [anchor=north west][inner sep=0.75pt]  [font=\normalsize]  {$t_{\gamma }$};
	\draw (302.62,98) node [anchor=north west][inner sep=0.75pt]  [font=\normalsize]  {$\arccos(t_{\gamma } \cdot n_{1});\, t_{\gamma } \cdot n_{1}=\mathcal{O}(h)$};

\end{tikzpicture}\vspace*{-0.15cm}

\caption{\hspace*{-0.1mm}Sketch \hspace*{-0.1mm}of \hspace*{-0.1mm}the \hspace*{-0.1mm}construction \hspace*{-0.1mm}as \hspace*{-0.1mm}described \hspace*{-0.1mm}in \hspace*{-0.1mm}Remark~\hspace*{-0.1mm}\ref{rem:soeren_trick_special_cases}~\hspace*{-0.1mm}with~\hspace*{-0.1mm}a~\hspace*{-0.1mm}discontinu-ity \hspace*{-0.1mm}set \hspace*{-0.1mm}$J_{z_T}$ \hspace*{-0.1mm}intersecting \hspace*{-0.1mm}an \hspace*{-0.1mm}element \hspace*{-0.1mm}$T$.
Part \textit{($\alpha$)} depicts 
 the~setting~of~Remark~\ref{rem:soeren_trick_special_cases}, while part \textit{($\beta$)} and part \textit{($\gamma$)} illustrate the cases \textit{(ii.a)} and \textit{(ii.b)} in Remark \ref{rem:soeren_trick_special_cases}.}
\label{fig:triangles}
	\end{figure}
	\end{remark}\newpage

	\section{Numerical experiments}\label{sec:experiments}
	
	\qquad In this section, we verify the theoretical findings of Section \ref{sec:non_sobolev} via~\mbox{numerical} experiments. To compare approximations to an exact solution,~we~impose~\mbox{Dirichlet} boundary \hspace*{-0.1mm}conditions \hspace*{-0.1mm}on \hspace*{-0.1mm}$\Gamma_D\hspace*{-0.2em}=\hspace*{-0.2em}\partial\Omega$, \hspace*{-0.1mm}though \hspace*{-0.1mm}an \hspace*{-0.1mm}existence \hspace*{-0.1mm}theory \hspace*{-0.1mm}is~\hspace*{-0.1mm}difficult~\hspace*{-0.1mm}to~\hspace*{-0.1mm}\mbox{establish}, in general. However, the error estimates derived in  Section \ref{sec:non_sobolev} carry over verbatimly with $\Gamma_N=\emptyset$ provided that a minimizer exists. 
	
	All \hspace*{-0.1mm}experiments \hspace*{-0.1mm}were \hspace*{-0.1mm}conducted \hspace*{-0.1mm}using \hspace*{-0.1mm}the \hspace*{-0.1mm}finite \hspace*{-0.1mm}element \hspace*{-0.1mm}software~\hspace*{-0.1mm}\textsf{FEniCS},~\hspace*{-0.1mm}cf.~\hspace*{-0.1mm}\cite{LW10}. All graphics are~generated~using the \textsf{Matplotlib} library, cf. \cite{Hun07}.
	
	\subsection{Experimental convergence rates}\label{sec:experiments_rates}
	\qquad All computations are based on  using the regularized discrete ROF~\mbox{functional}, i.e., for  $\varepsilon>0$ and $g\in L^2(\Omega)$, the functional $I^\varepsilon_h:\mathcal{S}^{1,cr}_D(\mathcal{T}_h)\to \mathbb{R}$, defined by
	\begin{align}
		I^\varepsilon_h(v_h):=\|\vert\nabla v_h\vert_\varepsilon\|_{L^1(\Omega)}+\frac{\alpha}{2}\|\Pi_h (v_h-g)\|_{L^2(\Omega)}^2\label{reg-rof}
	\end{align}
	for all $v_h\in \mathcal{S}^{1,cr}_D(\mathcal{T}_h)$, where $\vert \cdot\vert_\varepsilon\in C^1(\mathbb{R}^d)$ is the regularized modulus,~defined~by $\vert a\vert_\varepsilon:=(\vert a\vert ^2 +\varepsilon^2)^{\frac{1}{2}}$ for all $a \in \mathbb{R}^d$ and $\varepsilon>0$. On the basis of $0\leq \vert a\vert_\varepsilon- \vert a\vert\leq \varepsilon$ for  all $a\hspace*{-0.1em} \in\hspace*{-0.1em} \mathbb{R}^d$ and $\varepsilon\hspace*{-0.1em}>\hspace*{-0.1em}0$, for~the~minima~${u_h,u_h^\varepsilon\hspace*{-0.1em}\in\hspace*{-0.1em} \mathcal{S}^{1,cr}_D(\mathcal{T}_h)}$~of~${I_h,I_h^\varepsilon\hspace*{-0.1em}:\hspace*{-0.1em}\mathcal{S}^{1,cr}_D(\mathcal{T}_h)\hspace*{-0.1em}\to\hspace*{-0.1em} \mathbb{R}}$,~resp.,
	there holds
	\begin{align*}
		\frac{\alpha}{2}\|\Pi_h (u_h-u_h^\varepsilon)\|_{L^2(\Omega)}^2\leq \varepsilon\vert \Omega\vert.
	\end{align*}
	Thus, in order to bound the error $\|u-\Pi_hu_h\|_{L^2(\Omega)} $, it suffices to determine~the~\mbox{error} $\|u-\Pi_hu_h^\varepsilon\|_{L^2(\Omega)}$, e.g.,  for $\varepsilon=h$. The  iterative minimization of $\smash{I^h_h:\mathcal{S}^{1,cr}_D(\mathcal{T}_h)\to \mathbb{R}}$, i.e., for  $\varepsilon=h$, is realized using the unconditionally strongly stable semi-implicit discretized $L^2$--gradient flow from \cite{BDR18}, see also \cite[Section 5]{Bar20}.
	\begin{algorithm}[Semi-implicit discretized  $L^2$--gradient flow]\label{algorithm}
		Let $g_h\in \mathcal{L}^0(\mathcal{T}_h)$ and choose $\tau,\varepsilon_{stop}\hspace*{-0.1em} >\hspace*{-0.1em} 0$. Moreover, let
		$u^0_h\hspace*{-0.1em}\in\hspace*{-0.1em} \mathcal{S}^{1,cr}_D(\mathcal{T}_h)$ and set $k\hspace*{-0.1em}=\hspace*{-0.1em}1$.~Then,~for~${k\hspace*{-0.1em}\ge \hspace*{-0.1em} 1}$:
		\begin{description}[noitemsep,topsep=2pt,labelwidth=\widthof{\textit{(ii)}},leftmargin=!,font=\normalfont\itshape]
			\item[(i)] Compute $u_h^k\in \mathcal{S}^{1,cr}_D(\mathcal{T}_h)$ such that for every $v_h\in \mathcal{S}^{1,cr}_D(\mathcal{T}_h)$, there holds
			\begin{align*}
				\hspace*{-0.5cm}\big(d_tu_h^k,v_h\big)_{L^2(\Omega)}+\Bigg(\frac{\nabla_hu_h^k}{\big\vert \nabla_hu_h^{k-1}\big\vert_h},\nabla_hv_h \Bigg)_{\! \! L^2(\Omega;\mathbb{R}^d)}\!+\alpha \big(\Pi_hu_h^k-g_h,\Pi_hv_h\big)_{L^2(\Omega)}=0,
			\end{align*}
			where $d_tu_h^k:=\frac{1}{\tau}(u_h^k-u_h^{k-1})$ denotes the backward difference quotient.
		\item[(ii)]  Stop if $\big\|d_tu_h^k\big\|_{L^2(\Omega)}\!\leq \!\varepsilon_{stop}$; otherwise, increase $k\!\to\! k+1$~and~continue~with~\textit{(i)}.
		\end{description}
	\end{algorithm}
	
	It is shown in \cite[Proposition 3.4]{BDR18} and \cite[Proposition 5.3]{Bar20}, that Algorithm~\ref{algorithm} is unconditionally strongly stable, energy decreasing as well as converging, i.e., stops after finitely many iteration steps. To be more specific, for arbitrary $l\in \mathbb{N}$, one has the discrete energy estimate
	\begin{align*}
		I^h_h\big(u_h^l\big)+\tau\sum_{k=1}^l{\big\|d_tu_h^k\big\|_{L^2(\Omega)}^2}+\frac{\tau^2}{2}\sum_{k=1}^{l}{\int_{\Omega}{\frac{\big\vert d_t\nabla u_h^k\big\vert^2+\big(d_t\big\vert \nabla u_h^k\big\vert_h\big)^2}{\big\vert \nabla u_h^{k-1}\big\vert_h }\,\textup{d}x}}\leq I^h_h\big(u_h^0\big),
	\end{align*}
	which mainly results from $d_t\big\vert \nabla u_h^k\big\vert_h=\frac{1}{2}\big\vert \nabla u_h^{k-1}\big\vert_h^{-1}\big(d_t\big\vert\nabla  u_h^k\big\vert^2-\big(d_t\big\vert \nabla u_h^k\big\vert_h\big)^2\big)$.
	
	We will always employ the $h$--independent step-size $\tau = 1$ but the $h$--dependent \mbox{stopping} criteria $\big\|d_tu_h^k\big\|_{L^2(\Omega)}\leq\varepsilon_{stop}^h:=\frac{h}{20}$, i.e., $\big\|u_h^k-u_h^{k-1}\big\|_{L^2(\Omega)}\leq \frac{h}{20}$ as $\tau=1$.\newpage

	\begin{example}[Two disks problem]\label{two_disk_problem}
		Let $\Omega=(-1,1)^2\subseteq\mathbb{R}^2$, $r=0.4$, $\alpha =10$, and  $\tilde{g} := g\circ \Phi\in BV(\Omega)\cap L^\infty(\Omega)$, where $g:=\chi_{B_r^2(re_1)}-\chi_{B_r^2(-re_1)}\in BV(\Omega)\cap L^\infty(\Omega)$ and for some angle  $\phi\in [0,2\pi]$ and some vector $b_\gamma=(b_1,b_2)^\top\in \mathbb{R}^2$,
		\begin{align}
			\Phi(x):=\bigg[\begin{array}{c}
					\cos(\phi)(x_1-b_1)+\sin(\phi)(x_2-b_2)\\\cos(\phi)(x_2-b_2)-\sin(\phi)(x_1-b_1)
				\end{array}\bigg]\label{eq:two_disk_problem.1}
		\end{align}
		for all $x\hspace*{-0.1em}=\hspace*{-0.1em}(x_1,x_2)^\top\!\!\in\hspace*{-0.1em} \mathbb{R}^2$, i.e., $\Phi\hspace*{-0.1em}:\hspace*{-0.1em}\mathbb{R}^2\!\to\! \mathbb{R}^2$ 
		performs a rotation by $\phi$~and~a~shift~by~$b$. The same argumentation as in the proof of Proposition \ref{asym_primal_solution} demonstrates that the corresponding primal solution is given via $\tilde{u}:=u\circ \Phi =(1-\frac{2}{\alpha r})\tilde{g}\in BV(\Omega)\cap L^\infty(\Omega)$, where $u:=(1-\frac{2}{\alpha r})g\in BV(\Omega)\cap L^\infty(\Omega)$, cf.   Proposition \ref{asym_primal_solution}.

		For $z\in W^{\infty}(\textup{div};\Omega)$ defined as in Proposition \ref{asym_dual_solution}, we define the vector field $\tilde{z}:=\det(D\Phi)(D\Phi)^{-1}z\circ \Phi=(D\Phi)^{-1}z\circ \Phi\in W^{\infty}(\textup{div};\Omega)$. Then, resorting~to~properties of the contra-variant  Piola transform,~cf.~\cite[(2.1.71)]{BBF13},~we~find~that
		\begin{align}
			\textup{div}(\tilde{z})=\det(D\Phi)\textup{div}(z)\circ \Phi=\textup{div}(z)\circ \Phi=\alpha(\tilde{u}-\tilde{g})\quad\text{ in }\Omega.\label{eq:two_disk_problem.2}
		\end{align}
		We define the decomposition $\Omega^+_\Phi:=\Omega\cap \Phi(\mathbb{R}_{>0}\times \mathbb{R})$ and $\Omega^-_\Phi:=\Omega\cap \Phi(\mathbb{R}_{<0}\times \mathbb{R})$. Then,
		using that $\tilde{u}=0$ continuously on $b_\gamma+\mathbb{R}t_\gamma\cap \Omega$,~where~${t_\gamma=(-\sin(\phi),\cos(\phi))^\top}\!$,  and the transformation theorem, we further obtain that
		\begin{align}
			\begin{aligned}
			\vert D\tilde{u}\vert(\Omega)&=\vert D\tilde{u}\vert(\Omega^+_\Phi)+\vert D\tilde{u}\vert(\Omega^-_\Phi)
			\\&=	\vert Du\vert(\Omega^+)+\vert Du\vert(\Omega^-)
				\\&=	(u,\textup{div}(z))_{L^2(\Omega^+)}+(u,\textup{div}(z))_{L^2(\Omega^-)}
					\\&=	(u\circ \Phi,\textup{div}(z)\circ \Phi)_{L^2(\Phi^{-1}(\Omega))}=	(\tilde{u},\textup{div}(\tilde{z}))_{L^2(\Omega)},\end{aligned}\label{eq:two_disk_problem.3}
		\end{align}
		where we used in the last equality sign that $\textup{supp}(\tilde{u})\subseteq \Phi^{-1}(\Omega)\cap \Omega$.~Consequently,  if we combine \eqref{eq:two_disk_problem.2} and \eqref{eq:two_disk_problem.3} and  refer to the optimality conditions \eqref{optimality-relations},  then we find that
		 $\tilde{z}\in W^{\infty}(\textup{div};\Omega)$ is a dual solution to ${\tilde{u}\in BV(\Omega)\cap L^\infty(\Omega)}$. \mbox{Apparently},  $\tilde{z}\in W^{\infty}(\textup{div};\Omega)$ is  piece-wise Lipschitz continuous in the sense~of~Remark~\ref{rem:soeren_trick_special_cases} and its jump set is given via $J_{\tilde{z}}=b_\gamma+\mathbb{R} t_\gamma$. As a consequence, if for every $T\in \mathcal{T}_h$, either of the cases \textit{(ii.a)}--\textit{(ii.d)} in Remark \ref{rem:soeren_trick_special_cases} is satisfied, then the quasi-optimal rate $\smash{\mathcal{O}(h^{\frac{1}{2}})}$ is  guaranteed by Remark \ref{rem:soeren_trick_special_cases}, Lemma \ref{soeren_trick},~Lemma~\ref{dual_quasi_interpolant}~and~Theorem~\ref{non_Sobolev_error_estimate}.
	\end{example}

		\begin{example}[Four disks problem]\label{four_disk_problem}
		Let $\Omega=(-1,1)^2\subseteq\mathbb{R}^2$,~${r=0.4}$,~${\alpha =10}$,~and  
		\begin{align*}
			g:=\chi_{B_r^2(r,r)}+\chi_{B_r^2(-r,-r)}-\chi_{B_r^2(r,-r)}-\chi_{B_r^2(-r,r)}\in BV(\Omega)\cap L^\infty(\Omega).
		\end{align*}
		The same argumentation as for the proof of Proposition \ref{asym_primal_solution} shows that a minimum of \eqref{ROF-primal} is given via ${u:=(1-\frac{2}{r\alpha})g\in  BV(\Omega)\cap L^\infty(\Omega)}$.~A~\mbox{straightforward}  adaption of the proof of Corollary \ref{asym_primal_solution_corollary} implies that~any~dual~solution~${z \in W^\infty(\textup{div};\Omega)}$ is not $\theta$--Hölder continuous~at~${x=\pm r e_1}$~and~${x=\pm r e_2}$~if~${\theta>\frac{1}{2}}$.~Apart~from~that, arguing as in the proof of Proposition~\ref{asym_dual_solution}, we find that~an~example~of~a~dual~solution $\overline{z}\in W^\infty(\textup{div};\Omega)$ is~given~via 
		$\overline{z}(x):=
			\pm z(x\mp re_2)$ if $\pm x_2\ge 0$  
		for~all~$x=(x_1,x_2)^\top\in \Omega$, where $z\in W^\infty(\textup{div};\Omega)$ is defined as  in Proposition~\ref{asym_dual_solution}.~In~\mbox{addition}, if $\Phi\hspace*{-0.1em}:\hspace*{-0.1em}\mathbb{R}^2\hspace*{-0.1em}\to\hspace*{-0.1em} \mathbb{R}^2$ is defined as in Example \ref{two_disk_problem}, then for ${\tilde{g}\hspace*{-0.1em}:=\hspace*{-0.1em}g\circ \Phi\hspace*{-0.1em}\in\hspace*{-0.1em} BV(\Omega)\hspace*{-0.1em}\cap\hspace*{-0.1em} L^\infty(\Omega)}$, the primal solution is given via 
		${\tilde{u}:=u\circ \Phi\in  BV(\Omega)\cap L^\infty(\Omega)}$ and  a dual solution is given via  
		$\tilde{z}:=(D\Phi)^{-1}\overline{z}\circ \Phi\in W^\infty(\textup{div};\Omega)$. Apparently, $\tilde{z}\in W^\infty(\textup{div};\Omega)$ is piece-wise Lipschitz continuous in the sense~of~Remark~\ref{rem:soeren_trick_special_cases} and its jump~set~is~gi-ven via $J_{\tilde{z}}\!=\! b_\gamma+\mathbb{R}t_\gamma+\mathbb{R}n_\gamma$, where ${t_\gamma\!=\!(\!-\!\sin(\phi),\cos(\phi))^\top\!\!}$~and~${n_\gamma\!=\!(\cos(\phi),\sin(\phi))^\top\!\!}$. As a consequence, if for every $T\!\in\! \mathcal{T}_h$, either of the cases \textit{(ii.a)}--\textit{(ii.d)}~in~\mbox{Remark}~\ref{rem:soeren_trick_special_cases} is satisfied, then the quasi-optimal rate $\smash{\mathcal{O}(h^{\frac{1}{2}})}$ is guaranteed by Remark~\ref{rem:soeren_trick_special_cases}, Lemma~\ref{soeren_trick},~Lemma~\ref{dual_quasi_interpolant}~and~Theorem~\ref{non_Sobolev_error_estimate}.\vspace*{-4mm}
	\end{example}

\begin{figure}[h]
		\centering \vspace*{-0.35cm}
		
		\hspace*{-0.15cm}\includegraphics[width=12.5cm]{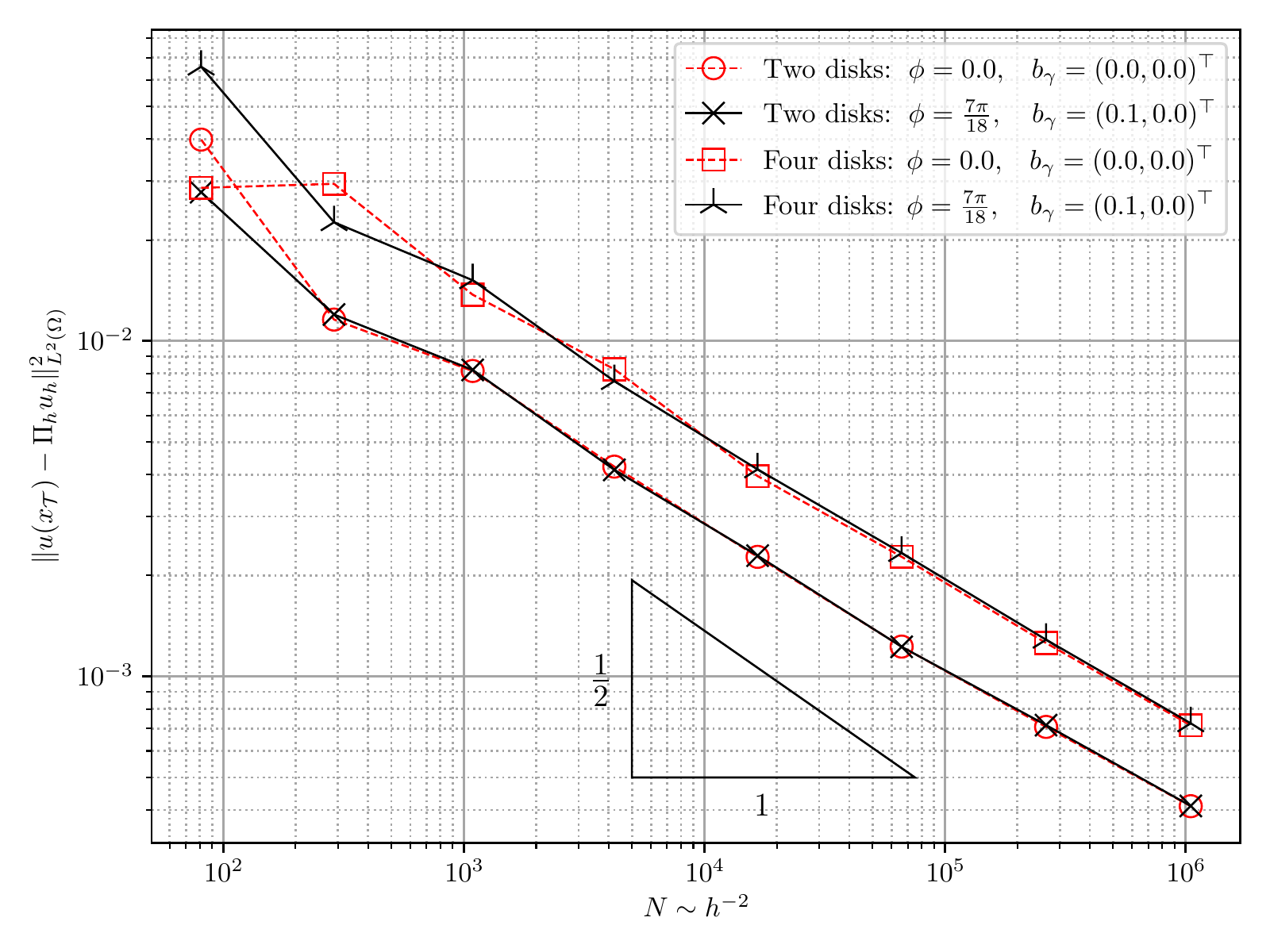}\vspace*{-0.75cm}
		
		\caption{Logarithmic plots for the experimental convergence rates of the error quantities \eqref{eq:L2-errors} 
			 in Example \ref{two_disk_problem} and in Example \ref{four_disk_problem}. The rate $\smash{\mathcal{O}(h^{\frac{1}{2}})}$ is observed.\vspace*{-12.5mm}}
		\label{fig:2Drates}
	\end{figure}

	The experimental convergence rates  in Figure \ref{fig:2Drates} are obtained~on~$k$--times~red-refined triangulations $\mathcal{T}_{h_k}$, $k=1,\dots,10$, of an initial triangulation $\mathcal{T}_{h_0}$ with two elements, i.e., $h_k=h_0 2^{-k}$ and $\varepsilon_{stop}^{h_k}=\frac{h_k}{20}$ for every ${k=1,\dots,10}$.~In~addition, for a simple implementation, we employ $\tilde{g}_{h_k}\in \mathcal{L}^0(\mathcal{T}_{h_k})$, defined by $\smash{\tilde{g}_{h_k}:=\tilde{g}(x_{\mathcal{T}_{h_k}})}$, where $x_{\mathcal{T}_{h_k}}|_T:=x_T$ for all $T\in \mathcal{T}_{h_k}$,
	instead~of~${g_{h_k}:=\Pi_{h_k}\tilde{g}\in \mathcal{L}^0(\mathcal{T}_{h_k})}$. However, since for each input data $\tilde{g}\in BV(\Omega)\cap L^\infty(\Omega)$ considered~in~this~section,~it~holds  $\|\tilde{g}-\tilde{g}_{h_k}\|_{L^1(\Omega)}\!\leq\! c h_k\vert \partial B_r^2(0)\vert $, the error estimate remains~valid.~To~be~more~\mbox{specific}, Figure \ref{fig:2Drates} contains logarithmic plots for the experimental convergence rates of the error quantities\vspace*{-1mm}
	\begin{align}
			\smash{\big\|u(x_{\mathcal{T}_{h_k}})-\Pi_{h_k}u_{h_k}\big\|_{L^2(\Omega)}^2,\quad k=3,\dots,10,}\label{eq:L2-errors}
	\end{align}
	versus the total number of vertices $N_k=(2^k+1)^2\sim h_k^{-2}$ for ${k=3,\dots, 10}$.~In~it,~we find  that the $L^2$--errors \eqref{eq:L2-errors} converge at the quasi-optimal convergence~rate~$\smash{\mathcal{O}(\hspace*{-0.05em}h^{\frac{1}{2}}\hspace*{-0.05em})}$. This behavior  is reported for both examples, i.e., Example \ref{two_disk_problem} and Example \ref{four_disk_problem}, for $\phi\hspace*{-0.1em}=\hspace*{-0.1em}0.0$ and $b_\gamma\hspace*{-0.1em}=\hspace*{-0.1em}(0.0,0.0)^\top$ as well as for  $\phi\hspace*{-0.1em}=\hspace*{-0.1em}\frac{7\pi}{18}$~and~${b_\gamma\hspace*{-0.1em}=\hspace*{-0.1em}(0.1,0.0)^\top}$.~Recall~that for $\phi= 0.0$ and $b_\gamma = (0.0,0.0)^\top$, the quasi-optimal rate~$\smash{\mathcal{O}(h^{\frac{1}{2}})}$~is~analytically guaranteed in both examples (cf. Example \ref{two_disk_problem} and Example \ref{four_disk_problem}).~Apart~from~that, we also could report the quasi-optimal rate $\smash{\mathcal{O}(\hspace*{-0.05em}h^{\frac{1}{2}}\hspace*{-0.05em})}$  for $d=3$, a uniform triangulation of $\Omega=(-1,1)^3$ and $g\in L^\infty(\Omega)\cap BV(\Omega)$ given via two or four touching~balls, with several rotations and shifts, for which no Lipschitz~continuous~dual~solution~exists.

	In Figure \ref{fig:four_disk_problem_plot},~the~\mbox{numerical} solution $u_{h_5}\hspace*{-0.1em}\in\hspace*{-0.1em} \mathcal{S}^{1,cr}(\mathcal{T}_{h_5})$ obtained in Example~\ref{four_disk_problem} and~\hspace*{-0.1mm}its~\hspace*{-0.1mm}\mbox{$L^2$--projection} \hspace*{-0.1mm}$\Pi_{h_5}u_{h_5}\!\in\! \mathcal{L}^0(\mathcal{T}_{h_5})$  \hspace*{-0.1mm}are \hspace*{-0.1mm}displayed \hspace*{-0.1mm}for~\hspace*{-0.1mm}${\phi\!=\!0.0}$~\hspace*{-0.1mm}and~\hspace*{-0.1mm}${b_\gamma\!=\!(0.0,0.0)^\top\!}$.
	Large gradients occur near the contact~points~of~the~disks, the midpoint values do not, however, show artifacts. 
	In Figure \ref{fig:discrete_dual_solution_plot},~the $L^2$--projection $\Pi_{h_5}z_{h_5}\in  \mathcal{L}^0(\mathcal{T}_{h_5})^2$ of the discrete dual solution ${z_{h_5}\!:=\!\nabla_{h_5} \!u_{h_5}\vert \nabla_{h_5} \!u_{h_5}\vert_{h_5}^{-1}\!+\!\frac{\alpha}{2}\Pi_{h_5}(u_{h_5}\!\!-\!g)(\textup{id}_{\mathbb{R}^2}\!-\!\Pi_{h_5}\textup{id}_{\mathbb{R}^2})}$ $\in\! \mathcal{R}T^0(\mathcal{T}_{h_5})$ \hspace*{-0.1mm}with \hspace*{-0.1mm}respect \hspace*{-0.1mm}to \hspace*{-0.1mm}the \hspace*{-0.1mm}regularized \hspace*{-0.1mm}ROF \hspace*{-0.1mm}functional~\hspace*{-0.1mm}\eqref{reg-rof}~\hspace*{-0.1mm}(cf.~\hspace*{-0.1mm}\mbox{\cite[\hspace*{-0.2mm}Section~\hspace*{-0.2mm}5]{BTW21}}) is displayed for ${\phi=0.0}$ and ${b_\gamma=(0.0,0.0)^\top}$.
	
\begin{samepage}
	\begin{figure}[h!]
	\hspace*{-1.5cm}\includegraphics[width=14.5cm]{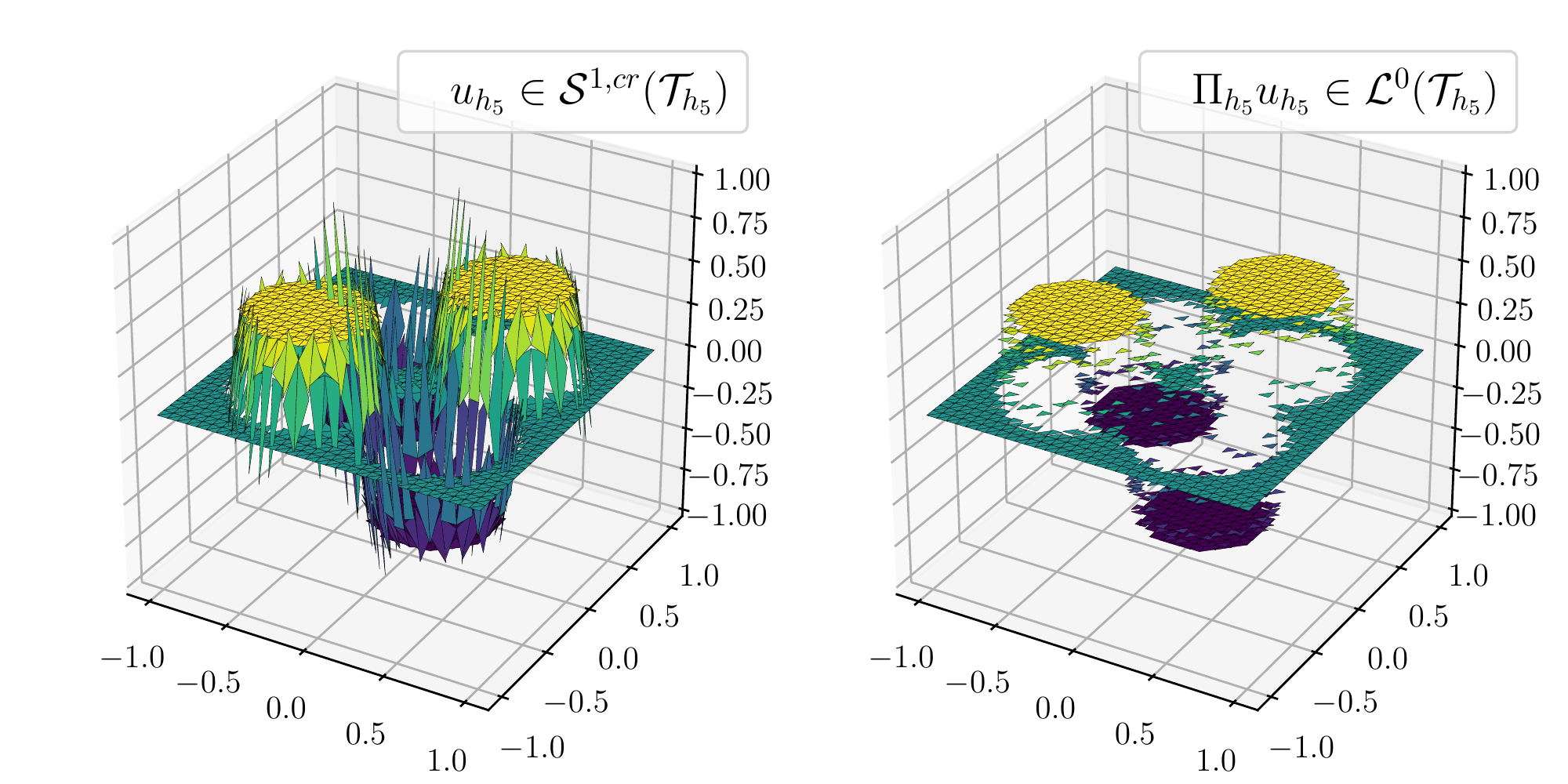}\vspace*{-0.1cm}
	\caption{Numerical solution $u_{h_5}\in \mathcal{S}^{1,cr}(\mathcal{T}_{h_5})$ in  Example \ref{four_disk_problem} displayed as piece-wise affine function (left) and via its $L^2$--projection~${\Pi_{h_5}u_{h_5}\in \mathcal{L}^0(\mathcal{T}_{h_5})}$~(right) for $r=0.4$, $\alpha =10$, $\phi=0.0$ and $b_\gamma=(0.0,0.0)^\top$. Large discrete~gradients~occur~near $\pm re_1$ and $\pm re_2$,  where no dual solution is $\theta$--Hölder continuous for ${\theta>\frac{1}{2}}$.\vspace*{-5mm}}
	\label{fig:four_disk_problem_plot}
\end{figure}
	\begin{figure}[h!]
		\hspace*{-1.65cm}\includegraphics[width=14.5cm]{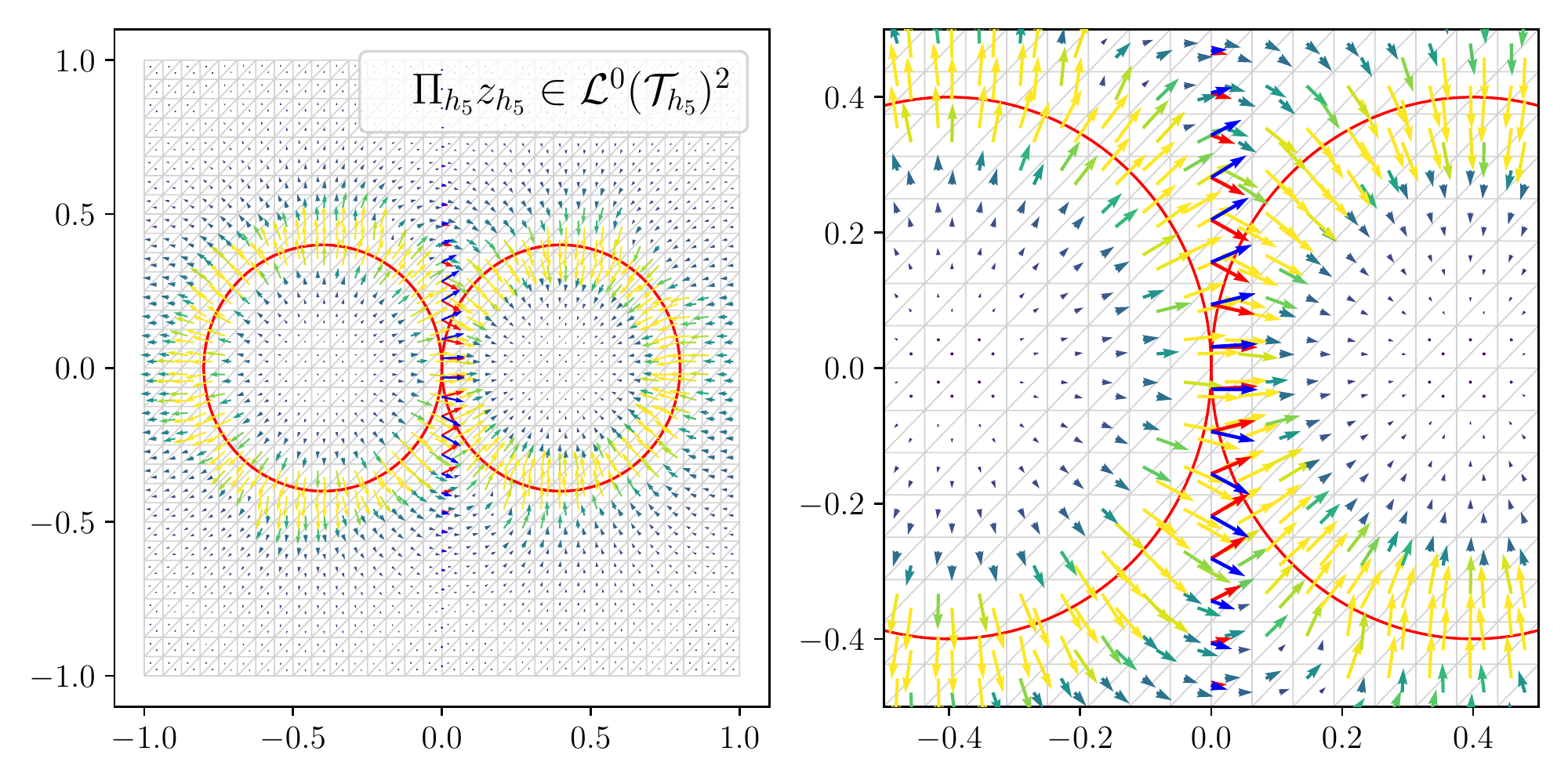}\vspace*{-0.1cm}
		\caption{$L^2$--projection $\Pi_{h_5}z_{h_5}\in  \mathcal{L}^0(\mathcal{T}_{h_5})^2$ of the discrete~dual~solution $z_{h_5}\in \mathcal{R}T^0(\mathcal{T}_{h_5})$ with  respect to the regularized ROF functional \eqref{reg-rof} (cf.~\mbox{\cite[Section~5]{BTW21}})  displayed for ${\phi=0.0}$ and ${b_\gamma=(0.0,0.0)^\top}$. The red and blue arrows represent the values of $z_{h_5}\in \mathcal{R}T^0(\mathcal{T}_{h_5})$ at the midpoints of element sides along the $\mathbb{R}e_2$--axis, i.e., $\lim_{\varepsilon\to 0} {z_{h_5}(x_S-\varepsilon e_1)}$ (blue arrows) and $\lim_{\varepsilon\to 0} {z_{h_5}(x_S+\varepsilon e_1)}$  (red arrows). Here, the different orientations of the arrows  indicate that $z_{h_5}\in \mathcal{R}T^0(\mathcal{T}_{h_5})$ approximates a discontinuous vector field -- empirically $z\in W^\infty(\textup{div};\Omega)$ defined in Proposition~\ref{asym_dual_solution}. Moreover, the red circles display the discontinuity set $J_u$ of the minimizer $u\in BV(\Omega)\cap L^\infty(\Omega)$ defined in Proposition \ref{asym_primal_solution}.\vspace*{-10mm}}
		\label{fig:discrete_dual_solution_plot}
\end{figure}\end{samepage}\newpage

\subsection{Experimental verification of condition \eqref{eq:soeren_trick}}

\begin{figure}[h]\vspace*{-0.25cm}
	\centering
	\hspace*{-0.25cm}\includegraphics[width=12.75cm]{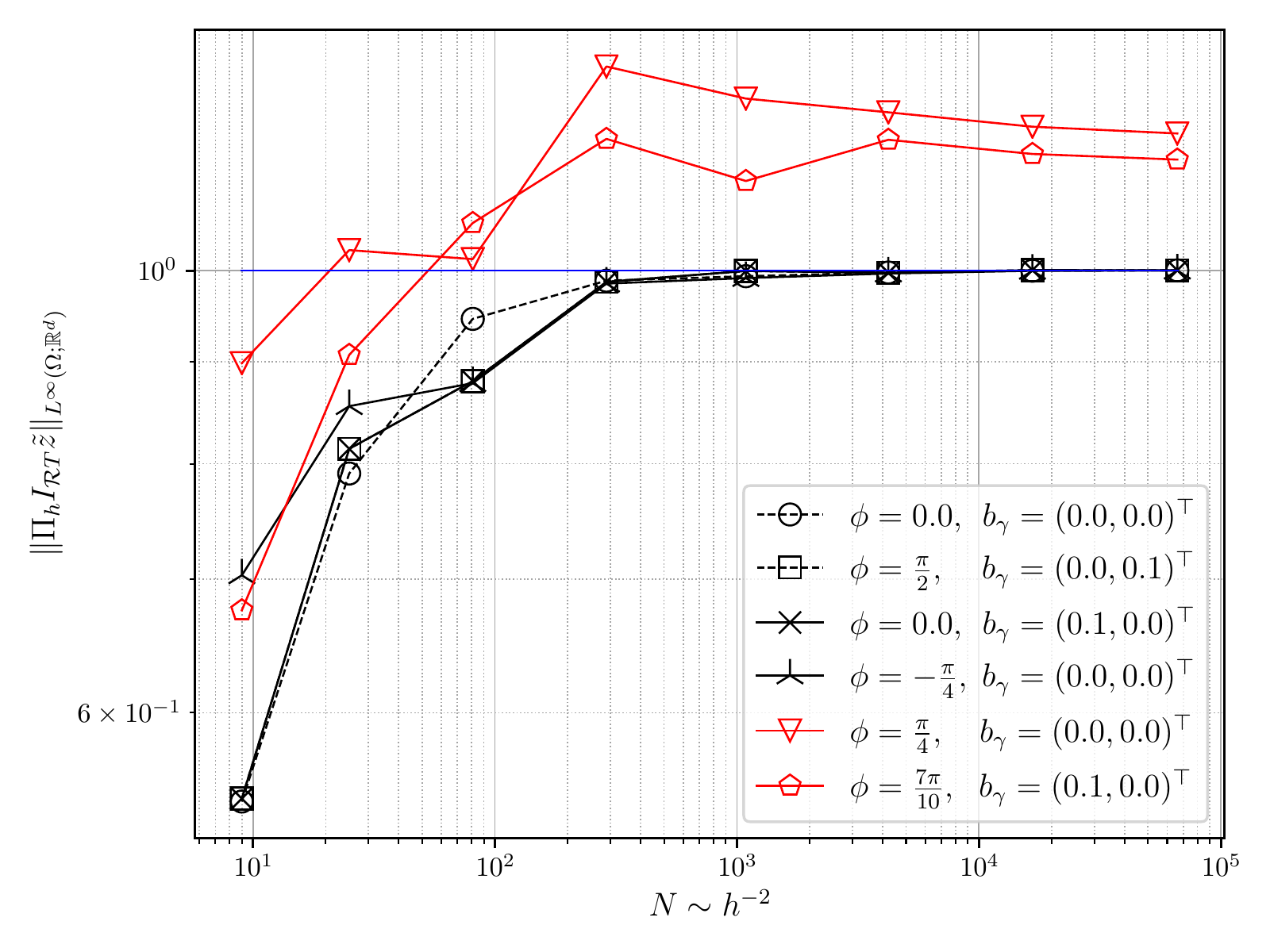}\vspace*{-0.5cm}
	\caption{Logarithmic plots of the quantities \eqref{eq:quantities1} 
		in Example \ref{two_disk_problem}.}
	\label{fig:interpolant}
\end{figure}

\qquad In this section, we examine whether the dual solutions given in Example~\ref{two_disk_problem} for every ${\phi \in [0,2\pi]}$ and $b_\gamma \in \mathbb{R}^2$ comply~with~condition~\eqref{eq:soeren_trick}~in~Lemma~\ref{soeren_trick},~which, in view of Lemma \ref{dual_quasi_interpolant} and Theorem \ref{non_Sobolev_error_estimate} yields a guarantee for the quasi-optimal convergence rate $\smash{\mathcal{O}(h^{\frac{1}{2}})}$. If we compute the quantities
\begin{align}
\|\Pi_{h_k}I_{\mathcal{R}T}\tilde{z}\|_{L^\infty(\Omega;\mathbb{R}^d)},\quad k=1,\dots,8,\label{eq:quantities1}
\end{align}
where $\tilde{z}\in W^\infty(\textup{div};\Omega)$ is defined as in Example \ref{two_disk_problem},
then we find that for $\phi =0.0$ and $b_\gamma =(0.0,0.0)^\top$, $\phi =\frac{\pi}{2}$ and  $b_\gamma =(0.0,0.1)^\top$, $\phi =0.0$ and  $b_\gamma =(0.1,0.0)^\top$,  and  $\phi =-\frac{\pi}{4}$ and $b_\gamma =(0.0,0.0)^\top$, there exists
a constant $c_z>0$ -- presumably,  one has that $c_z=1$ -- such that  for $k=1,\dots,8$, there holds
\begin{align}
	\|\Pi_{h_k}I_{\mathcal{R}T}\tilde{z}\|_{L^\infty(\Omega;\mathbb{R}^d)}\leq 1+c_zh_k.\label{estimate}
\end{align}
These results confirm the findings in Remark \ref{rem:soeren_trick_special_cases}  as they fall within one of the cases \textit{(ii.a)}--\textit{(ii.d)}. Apart from that, for $\phi =\frac{\pi}{4}$ and $b_\gamma =(0.0,0.0)^\top$ as well as for $\phi =\frac{7\pi}{18}$ and $b_\gamma =(0.0,0.0)^\top$, we cannot report the existence of a constant $c_z>0$ such that \eqref{estimate}   holds.  This behavior can also be easily~\mbox{predicted} analytically by resorting to the formula \eqref{eq:soeren_trick_special_cases.3}.
All results can be found in~\mbox{Figure}~\ref{fig:interpolant}, which displays  the quantities \eqref{eq:quantities1} 
	versus the total number of vertices ${N_k\!=\!(2k+1)^2\!\sim\! h_k^{-2}}$~for~${k\!= \! 1,\dots,8}$.

	\vspace*{-0.15cm}

Possible explanations for the observed quasi-optimal rate $\smash{\mathcal{O}(h^{\frac{1}{2}})}$ for $\phi =\frac{\pi}{4}$ and $b_\gamma =(0.0,0.0)^\top$ as well as for $\phi =\frac{7\pi}{18}$ and $b_\gamma =(0.0,0.0)^\top$, even though \eqref{estimate} could not be reported, might be that this violation is merely pre-asymptotic or
 occurs only along the interface $(b_\gamma+\mathbb{R}t_\gamma)\cap \Omega$ (the latter, we observed experimentally), that the proofs presented are still sub-optimal, or that there exists an alternative dual solution for which \eqref{estimate} can be reported.\vspace*{-6mm}\newpage

		\providecommand{\MR}[1]{}
		\providecommand{\bysame}{\leavevmode\hbox to3em{\hrulefill}\thinspace}
		\providecommand{\noopsort}[1]{}
		\providecommand{\mr}[1]{\href{http://www.ams.org/mathscinet-getitem?mr=#1}{MR~#1}}
		\providecommand{\zbl}[1]{\href{http://www.zentralblatt-math.org/zmath/en/search/?q=an:#1}{Zbl~#1}}
		\providecommand{\jfm}[1]{\href{http://www.emis.de/cgi-bin/JFM-item?#1}{JFM~#1}}
		\providecommand{\arxiv}[1]{\href{http://www.arxiv.org/abs/#1}{arXiv~#1}}
		\providecommand{\doi}[1]{\url{http://dx.doi.org/#1}}
		\providecommand{\MR}{\relax\ifhmode\unskip\space\fi MR }
		\providecommand{\MRhref}[2]{%
			\href{http://www.ams.org/mathscinet-getitem?mr=#1}{#2}
		}
		\providecommand{\href}[2]{#2}
		{\setlength{\bibsep}{0pt plus 0.0ex}
	}

\end{document}